\documentclass[a4paper,12pt]{amsart}
\pdfoutput=1
\usepackage[protrusion=true,expansion=true]{microtype}
\usepackage{xcolor}
\usepackage[T1]{fontenc} 
\usepackage{mathtools}
\usepackage[all]{xy}
\usepackage[english]{babel}
\usepackage{frcursive} 
\usepackage[applemac]{inputenc}
\usepackage{eucal}
\usepackage{graphicx}
\usepackage{epsfig}
\usepackage{mathrsfs}
\usepackage{amssymb}
\usepackage{amsxtra}
\usepackage{enumerate}
\input cyracc.def 
 
\newtheorem*{theorem}{Theorem}
\newtheorem*{theo}{Theorem}
\newtheorem*{proposition}{ Proposition}
\newtheorem*{lemma}{ Lemma}

\newtheorem*{corollary}{Corollary}

\newtheorem*{definition}{Definition}

\theoremstyle{remark}

\def \1{\mathbb {1}}
\def \RM{\mathbb {R}}%        corps des reels
\def \NM{\mathbb{N}}%        entiers naturels
\def \Nb{\overline{\mathbb{N}}}%
\def \ZM{\mathbb{Z}}%        entiers relatifs
\def \TM{\mathbb{T}}
\def \CM{\mathbb{C}}%        nombres complexes

\def \QM{\mathbb{Q}}%        nombres rationnels

\def \SM {\mathbb{S}} 

 \def \Hom {{\rm Hom}}

\def \ord {{\rm ord\,}}

\def \Der {{\rm Der\,}}

\def \p {{\rm exp\,}}
\def \Id {{\rm Id\,}}

\def \d{\partial}%derivee partielle
\def\dt{\delta} 
\def\a{\alpha}
\def\b{\beta}
\def\e{\varepsilon}  
\def\g{\gamma}

\def\L{\Lambda}
\def\p{\varphi}

\def\G{\Gamma}   
\def\D{\Delta}
\def \s{\sigma}

\def \to{\longrightarrow} 
\def \w{\wedge}

\def \< {{\langle }}
\def \> {{\rangle }}
\def \( {\left( }
\def \) {\right) }

\newcommand{\Bt}{{\mathcal B}}
\newcommand{\Ct}{{\mathcal C}}

\newcommand{\Ht}{{\mathcal H}}

\newcommand{\Lt}{{\mathcal L}}

\newcommand{\Ot}{{\mathcal O}}
\newcommand{\Pt}{{\mathcal P}}
\newcommand{\Qt}{{\mathcal Q}}

\newcommand{\St}{{\mathcal S}}

\newcommand{\Vt}{{\mathcal V}}

\newcommand{\Xt}{{\mathcal X}}

\newcommand{\Sform}{\bf S}
\title[Quasi-periodic motions on symplectic tori]{Quasi-periodic motions on symplectic tori}
\author[M. Garay A. Kessi D. van Straten N. Yousfi]{Mauricio GARAY, Arezki KESSI,\\ Duco van Straten and  Nesrine YOUSFI.}
\begin{document}

%\normalsize
%\authorfootnotes

\begin{abstract}
The KAM (Kolmogorov-Arnold-Moser) theorems guarantee the stability of 
quasi-periodic invariant tori by perturbation in some Hamiltonian systems.
Michel Herman proved a similar result for quasi-periodic motions for 
$k=2d-2$-dimensional involutive manifolds in Hamiltonian systems with $d$ 
degrees of freedom. In this paper, we give a similar result for quasi-periodic 
motions on symplectic tori, i.e. the case $k=2d$.
\end{abstract}
\maketitle

\section*{Introduction}

The theorems of Kolmogorov, Arnold and Moser guarantee, under appropriate
conditions, the stability of quasi-periodic motions on Lagrangian invariant
tori~\cite{Arnold_KAM,Fejoz_KAM,Kolmogorov_KAM,Moser_KAM,Poschel_lecture,Russmann_degeneracy,Sevryuk_KAM}. The phase space then is a symplectic manifold of dimension $2d$ and
these KAM-tori have dimension~$d$. One may ask about the stability of other 
types of quasi-periodic motions that fill out a torus of dimension $k \neq d$. 
The case of such movements on a lower dimensional torus ($k <d$) was already 
considered by Moser back in 1967 \cite{Moser_convergent}. But it was only 
in 1991 that Herman exhibited examples of Hamiltonian systems with stable 
quasi-periodic motion on a $k=2d-2$ dimensional co-isotropic torus~\cite{Herman_KAM,Yoccoz_Herman}. 

The aim of this paper to formulate and prove a stability theorem for motions on 
symplectic tori, which can be interpreted as the case $k=2d$.
We believe that the similar arguments may be used to handle the cases for any
$d <k\le 2d$, but to keep the paper as simple as possible we decided to concentrate on the extremal case of quasi-periodic motion on a symplectic torus. 
 
In the Lagrangian or isotropic case, variations of the symplectic
structure on the phase space are usually not considered, as they 
can be cancelled by an adequate change of variables.

Herman used the fact that the symplectic structure can be varied 
along general co-isotropic manifolds and hence implicitly considered 
variations of the symplectic structure.

As we shall see, the phenomenon discovered by Herman holds not only for 
Hamiltonian flows, but can be generalised to quasi-periodic motions
generated by more general symplectic vector fields. Although these 
are technically more difficult to handle, the statements obtained are 
not only more general, but also simpler to formulate.

We found that a purely na\"ive point of view and a direct attack on the 
problem leads into a morass of technical difficulties that has prevent us 
from giving a convincing presentation of a proof along these lines.
To obtain a transparent proof, we found it convenient to use a more
abstract framework and adopt the language of {\em Kolmogorov spaces} of 
abstract KAM theory developed in \cite{KAM_book}. These analytic methods are inspired from the work of Douady and Hauser~\cite{Douady_these,Hauser_these}. We sketch the main
features of this framework, but the paper will be essentially self-contained 
and full proofs are given, except for some simple lemmas.\\

The structure of the paper is as follows. In the first section we formulate
our main theorem on the stabilty of real analytic symplectic vector fields
on a symplectic torus under an arithmetic condition on the frequencies which is much weaker than the usual Diophantine condition.
As usual in a real analytic context, we will use holomorphic tools on
a complex analytic neighbourhood of the torus. This is explained in some detail 
in the second section and serves to fix some of the notations we use.
We introduce formal parameters that describe the {\em detuning} of the
frequencies and {\em perturbation} of the symplectic form. In the third section
we give an algebraic description of an almost quadratic iteration in 
terms of a certain ring $R$, that brings our vector field formally into 
normal form.
In section six this iteration scheme is lifted to an iteration scheme
in a functional analytic context, so that the issue of convergence can be 
adressed. In order to do so, we have to make some careful preparations. In 
section four we sketch out an abstract framework for handling families of 
Banach spaces parametrised by partially ordered sets that we call 
{\em Kolmogorov spaces}. These naturally arise in situations where one has
 to deal with functions and vector fields defined over domains that shrink
during an iteration proces, due to the appearance of small denominators. 
Most relevant is the general theorem formulated in 4.5.2, which is crucial 
in the later part of the paper to control the norm estimates, in particular
for the exponentials of vector fields that we use.  In section five
we describe the neighbourhoods of the resonance hyperplanes that
need to be removed in order to control the small denominators and we
introduce the precise Kolmogorov spaces of holomorphic functions we will use
in our proof. In section seven we lift the iteration scheme on the level 
of Kolmogorov  spaces. The required estimates now are all automatic and the proof is completed with ease, by comparison with a simple model iteration.\\

\section{Quasi-periodic movement on a symplectic torus}

\subsection{}
We will describe our basic set-up und formulate the main result of this
paper. By {\em quasi-periodic motion} on a torus
\[\TM^{n}:=(\RM/\ZM)^{n}\]
with coordinates $x_1,\dots,x_{n}$ we mean the flow of a 
constant vector field
\[ X:=\sum_{i=1}^{n} \nu_i \frac{\partial}{\partial x_i}.\] 
If the components of the {\em frequency vector}
\[ \nu:=(\nu_1,\nu_2,\ldots,\nu_n) \in \RM^n\]
are independent over $\QM$, the orbits of $X$ are dense in $\TM^n$
and we say $\nu$ is {\em non-resonant}.
We will refer to these constant vector fields as {\em quasi-periodic vector
fields}. It is clear that by an arbitrary small perturbation $S$ of the vector 
field $X$ we can create a vector field $X+S$ that is not conjugate to $X$,
so in this sense the dynamics generated by the vector fields $X$ is not stable.
But for a special class of perturbations, such a conjugation will be seen to
be possible.
\subsection{} 
Consider the case where $n=2d$ and equip the torus $\TM:=\TM^{2d}$ 
with a constant symplectic form
$$\omega=\sum_{1 \le i<j \leq 2d} \omega_{ij} dx_i \w dx_j,\ \omega_{ij} \in \RM.$$
A $C^{\infty}$ vector field $X$ on $\TM$ is called a 
{\em symplectic vector field} if the time $t$ flow $\Phi_t$ of 
$X$ preserves the symplectic form:
\[ \Phi^*_t(\omega)=\omega .\]
This is equivalent to the infinitesimal condition
\[ L_X(\omega)=0,\]
where 
\[ L_X=d\iota_X+\iota_X d\]
denotes the {\em Lie-derivative}. As $d\omega=0$, this is equivalent to the
statement that the one-form $\iota_X \omega$ symplectically dual to $X$ is
{\em closed}. In the particular case where the form is {\em exact}, the
corresponding vector field $X$ satifies
\[\iota_X\omega =dh\] 
for some function $h$ on $\TM$ and we say that the field is {\em Hamiltonian}, or more precisely that the vector field is 
{\em associated to the Hamiltonian function $h$.}  
Obviously, constant vector fields are symplectic, but not Hamiltonian.

We will now seek a statement expressing the {\em stability of quasi-periodic motions under perturbation with symplectic vector fields}. For this to work properly, we will also allow for perturbations of the chosen symplectic form.\\

\subsection{} The space $\Sform$ of constant symplectic forms on the torus $\TM$ can be identified with those on $\RM^{2d}$ and thus forms the open subset of non-degenerate skew-symmetric 2-forms
$$\Sform \subset \L^2(\RM^{2d})^* \simeq \RM^{d(2d-1)}.$$  
The fibre of the fibration
$$\pi: \TM \times \Sform \to \Sform $$
over the point
$$\omega=\sum_{i<j}\omega_{ij}dx_i \w dx_j \in \Sform $$ 
describes the torus $\TM$ with symplectic form $\omega$.

In this situation, the relative 1-forms are of the form
$$\sum_{i=1}^{2d} a_i(x,\omega)dx_i \in \Omega^1_\pi,$$
and dually, the relative vector fields are of the form
\[\sum_{i=1}^{2d} b_i(x,\omega)\d_{x_i} \in \Theta_{\pi} .\]

The interior product with the symplectic form $\omega$ induces an isomorphism
$$\Theta_{\pi} \to \Omega^1_\pi , \;\;\; X \mapsto \iota_X \omega :=\omega(X,-). $$

A {\em quasi-periodic vector field depending on parameters $\omega \in V \subset \Sform$}
is a vector field on $\TM \times V$ of the form:
$$X_\nu=\sum_{i=1}^{2d}\nu_i(\omega) \d_{x _i}.$$
where
$$\nu:=(\nu_1,\dots,\nu_{2d}): V \to \RM^{2d}  $$
is called {\em the frequency map} of the vector field. 
Such a vector field generates quasi-periodic movements on the fibres 
$\TM \times \{\omega\}$, for all $\omega \in V$, with frequencies that 
may depend on $\omega$. 

\subsection{}
Our aim is to prove the following {\em Stability Theorem:}

\begin{theo} Consider a quasi-periodic motion defined by a real analytic vector field
$$X_\nu=\sum_{i=1}^{d}\nu_i(\omega) \d_{x _i}$$ 
defined on a neighbourhood  $\TM \times V$ of a symplectic torus $\TM  \times \{ \omega^0 \}$. 
Assume that 
\begin{enumerate}
\item[(A)] the vector $\nu(0)$ satisfies a subquadratic arithmetic condition. 
\item[(B)] the map $$\nu: V \to \RM^{2d},\;\omega \mapsto (\nu_1(\omega),\dots,\nu_{2d}(\omega))$$ 
is a submersion.
\end{enumerate}

Then: For any real analytic symplectic vector field  $X$ sufficiently close to $X_\nu$, there exists a set $M \subset V$ of positive measure parametrising tori that carry a motion conjugate by a symplectomorphism to a quasi-periodic one.
\end{theo}

The topology on the  space of analytic vector fields will be reviewed in later 
sections of the paper. 

\subsection{}
We now spell-out the arithmetic condition of the theorem. By definition, a frequency vector $\nu=(\nu_1,\nu_2,\ldots,\nu_{2d})$ is called 
{\em non-resonant} if the scalar product
\[ (\nu,J):=\sum_{k=1}^{2d} \nu_k J_k,\;\;\; J=(J_1,J_2,\ldots,J_{2d}) \]
is non-zero for all $J \in \ZM^{2d}\setminus\{0\}$. But being non-zero,
this quantity can become arbitrary small for certain $J$, if we allow $|J|$ 
to become large. As during the iteration one has to divide
by this quantity, these {\em small denominators} have a dangerous effect on the
convergence and must be controlled.\\
A convenient way to quantify such small denominators is by the so-called 
{\em arithmetic sequence} 
$$\sigma(\nu)=(\sigma(\nu)_k)$$ 
attached to a vector $\nu \in \CM^{2d}$. This falling sequence is defined by setting
$$\s(\nu)_k :=\min \{ |(\nu,J)|: J \in \ZM^{2d} \setminus \{ 0 \}, \| J \| \leq 2^k \} .$$
If we collect all frequency vectors $\nu$ for which this arithmetic sequence is bounded from below by a given falling sequence $a=(a_k),$ we obtain what we call the {\em arithmetic class} of $a$,  defined as
\[\Ct(a):=\cap_{m=0}^{\infty} \Ct_m(a),\;\;\Ct_m(a)=:\{ \nu \;\; |\;\; \sigma(\nu)_k \ge a_k, k=1,2,\ldots,m\}.\]

The Cantor-like set $\Ct(a)$ could be called a {\em swiss cheese set}, as it is obtained be removing smaller and smaller neighbourhoods around the dense collection of hyperplanes 
$(\omega,J)=0$, $0 \neq J \in \ZM^{2d}$. Obviously, one has
\[ a' \le a,\;\;\;\Longrightarrow \Ct(a) \subset \Ct(a').\]
It is an elementary fact that for any $\nu \in \Ct(a)$, we can find $a' \le a$ such that $\Ct(a')$ has density $=1$ at the point $\nu$, see e.g. \cite{arithmetic}.\\

One says that $\nu$ is {\em Diophantine} if it satisfies a {\em Diophantine condition}, meaning that there exist constants $C$ and $N$ such that 
\[ |(\nu,J)| \ge \frac{C}{|J|^N} .\]
This means that $\nu \in \Ct(a)$ where $a$ is a falling geometrical 
sequence.
Diophantine conditions appear often in dynamical systems, but after the work of Bruno it appeared that in many case the condition can be relaxed~\cite{Bruno}. In our theorem we will need only much weaker condition than the Diophantine one: the {\em subquadratic arithmetic condition} that we explain now.
Note that $a_n=q^{\a^n}$ solves the iteration
\[   a_{n+1} =(a_n)^\a  \]
and this motives that following definition.
 
\begin{definition}
A strictly increasing sequence of positive numbers $a=(a_n)$ is called 
{\em positively subquadratic with exponent $\a$,} if $\a \in]1,2[$ and
there exists $A,B \in \RM_{>0}$ such that for all $n \in \NM$ one has
$$a_n \leq Ae^{B\a^n}. $$
We denote the set of such sequences by $\SM^+(\a)$.
Similarly, a strictly decreasing sequence of positive numbers $a=(a_n)$ is
called  {\em negatively subquadratic with exponent $\a$} if $\a \in ]1,2[$
and there exists $A,B \in \RM_{>0}$ such that for all $n \in \NM$ one has
$$a_n \geq Ae^{-B\a^n}.$$ 
We denote the set of such sequences by $\SM^-(\a)$.
\end{definition}

Clearly, one has $\SM^+(\a) \subset \SM^+(\b)$ if $\a \le \b$, so that the set 
\[\SM^+:= \bigcup_{\a \in ]1,2[} \SM^+(\a)\]
is filtered be the sets $\SM^+(\a)$. Note also that the product of
such subquadratic sequences is again subquadratic and one has
\[\SM^+(\a) \cdot \SM^+(\b) \subset \SM^+(\max(\a,\b))\]
Given $a\in \SM^+$, we call the {\em order of  $a$}, denoted $\ord(a)$ the infimum of its exponents.
Taking the multiplicative inverse  $(a_n) \mapsto (1/a_n)$ interchanges $\SM^+(\a)$ and $\SM^-(\a)$, so that one has corresponding properties.
\[\SM^-(\a) \cdot \SM^-(\b) \subset \SM^-(\max(\a,\b))\]

\begin{definition}
A frequency vector $\nu \in \CM^n$ is said to be {\em subquadratic} if $\s(\nu) \in \SM^-$. 
\end{definition}

For a subquadratic sequence the infinite product
\[ a_{\Pi}:=a_0 a_1^{1/2}a_2^{1/2^2}\ldots=\prod_{k=0}^{\infty} a_k^{1/2^k}\] 
converges to a strictly positive number or equivalently if
$$\sum_{k\geq 0} \left|\frac{\log a_k}{2^k}\right| <+\infty. $$
The sequences satisfying this last condition are called {\em Bruno sequences}. So subquadratic sequences form a subset of Bruno sequences. 
Clearly, if  $\nu$ is Diophantine, then  $\s(\nu)$ bounded by a gemetrical 
sequence, so these are in particular subquadratic sequences with exponent 
equal to $1$.
%%%%%%%%%%%%%%%%%%%%%%%%%%%%%%%%%%%%%%%%%%%%%%%%%%%%%%%%%%%%%%%%%%%%%%%%%%%%%%
\section{The analytic torus}
As usual in KAM theory, our theorem is proved using {\em normal form} 
techniques and more precisely we will use ideas from {\em parametrised KAM theory}~\cite{Broer_Huitema_Sevryuk_book,Broer_Huitema_Takens,Broer_Huitema_Sevryuk_families}
(see also \cite{Gonzalez_Haro_DelaLLave}). Our normal form is achieved using a specific iteration scheme. As usual in a real analytic context, 
we use complexification and construct the iteration in appropriate open subset 
in the complex domain. As all constructions can be done compatible with the
underlying real structure, nothing is lost and much is won by doing so.
Before we discuss the iteration itself, we set up the basic analytic notions on the torus.
\subsection{Analytic functions on the torus}

\subsubsection{} 
The exponential map 
\[ \RM^{n} \to (\CM^*)^{n},\;\;\;(x_1,x_2,\ldots,x_{n}) \mapsto (z_1,z_2,\ldots,z_{n})\]
with
\[ z_j=e^{ 2 \pi i x_j}, \;\;\;j=1,2,\ldots,n\]
defines an embedding of the real torus $\TM^{n}=(\RM/\ZM)^n$ in the algebraic torus $(\CM^*)^{n}$ as the product of unit circles $|z_i|=1$. We will identify $\TM^n$ with
this subset of $(\CM^*)^n$ and describe functions,
vector fields, differential forms on the torus using the complex coordinates
$z_j$. 
\subsubsection{}
The simplest functions on the torus $\TM^n$ are the {\em Laurent-polynomials},
which are of the form
\[ f=\sum_{I \in \ZM^{n}} a_Iz^I, \]
where only finitely many $a_I \neq 0$. Here and in the sequel we use 
the {\em multi-index notation} and write 
\[ z^I=z_1^{i_1} z_2^{i_2}\ldots z_n^{i_n},\;\;I=(i_1,i_2, \dots,i_n),\;\;\;|I|=|i_1|+|i_2|+\dots+|i_n|,\;\; \textup{etc}. \]
The set of all Laurent-polynomials forms a ring denoted by $\CM[z,z^{-1}]$;
these correspond precisely to trigonometric polynomials when written in $x$.

\subsubsection{}
The analytic functions on $\TM^n$ are identified with the ring 
\[ A:= \CM\{ z,z^{-1} \}\] 
of {\em analytic Fourier series} and consist of functions which are 
holomorphic on an open neighbourhood of the real torus $\TM^n$. A fundamental
system of such neighbourhoods is provided by the sets 
\[T_r:=\{ z \in \CM^n\;|\; e^{-r} < |z_i| <e^r\}, \;\;\;t \in ]0,\infty].\] 
So  $\CM\{z,z^{-1}\}$ is the union of Banach spaces $\Ot^b(T_r)$, the space 
of {\em bounded holomorphic functions} on $T_r$.
As a result, $A$ carries a natural Fr\'echet space structure. Elements
of $\CM\{z,z^{-1}\}$ are represented  by power series of the form
\[\sum_{I \in \ZM^{n}} a_Iz^I \text{ such that }  \sum_{I \in \ZM^{n}} |a_I|e^{r|I|}<+\infty\]
for some $r>0$. Alternatively, $\CM\{z,z^{-1}\}$ can be seen as the union
of Banach spaces $\Ot^k(T_r)$ of holomorphic functions on $T_r$ with $C^k$-
extension to $\overline{T_r}$, or of the Hilbert spaces $\Ot^h(T_r)$ of square 
integrable homolorphic functions on $T_r$.
\subsubsection{}
One has inclusions 
\[ \CM[z,z^{-1}] \subset \CM\{z,z^{-1}\} \subset \CM[[z,z^{-1}]], \]
where on the right hand side  we have the vector space of {\em formal Laurent series}
\[  F:=\sum_{I \in \ZM^n} a_I z^I,\]
where the coefficients $a_I$ are completely arbitrary. The Cauchy product of
two such formal Laurent series is usually not defined, so the ring structure of
$\CM\{z,z^{-1}\}$ does not extend to $\CM[[z,z^{-1}]]$.  However, the 
{\em Hadamard-product} $\star$ obtained by
coefficient-wise multiplication 
\[ \sum_{I \in ZM^n} a_I z^I \star \sum_{I \in \ZM^n} b_I z^I =\sum_{I\in \ZM^n} a_Ib_I z^I\]
defines a well-defined operation on formal Laurent series:
\[ \star: \CM[[z,z^{-1}]] \times \CM[[z,z^{-1}]] \to \CM[[z,z^{-1}]].\]

For a  formal Laurent series 
\[h =\sum_{I \in \ZM^n} h_I z^I \in \CM[[z,z^{-1}]]\] 
the operation $h \star :\CM[[z,z^{-1}]] \to \CM[[z,z^{-1}]]$ 
maps $\CM[z,z^{-1}]$ to itself, but in general does not preserve 
$\CM\{z,z^{-1}\}$
if the coefficient $h_I$ grow too fast for $|I| \to \infty$. 
However, if the coefficients satisfy an estimate of the form
\[ |h_I| \le C |I|^N,\]  
for some $C,N$, then $h \star$ maps $\CM\{z,z^{-1}\}$ to itself.
\subsubsection{}
For each subset $A \subset \ZM^n$ there is a canonical {\em truncation operator}
\[ [-]_A : \CM[[z,z^{-1}]] \to \CM[[z,z^{-1}]], f \mapsto [f]_A ,\]
where we keep only the the monomials of $f$ whose exponent appear in $A$:
\[[f]_A:= \sum_{I \in \ZM^{n}\cap A} a_Iz^I .\]
So we are dealing with a special case of the Hadamard product and so these truncation operators map
$\CM\{z,z^{-1}\}$ to itself.
If $A$ is a finite set, $[f]_A$ belongs to $\CM[z,z^{-1}]$.
For $j,k \in \NM$ we put
\[ [f]^k_j:=[f]_{A_{j,k}},\;\;\;A_{j,k}:=\{ I \in \ZM^{n}\;|\;\;\;j \leq \;|I| < k\}.\]
We also write $[f]_j$ in case $k$ is infinite and $[f]^k$ in case $j=0$.

\subsection{Analytic vector fields on the torus}

\subsubsection{}
By an analytic vector field on $\TM^n$ we mean a vector field
\[ X=\sum_{j=1}^{n} a_j(x) \d_{x_j} , \]
where the coefficients $a_j(x) \in A$, so are analytic functions on the torus.
Expanding these coefficients in Fourier series and using 
$$\frac{1}{2\pi i} \d_{x_j}= z_j\d_{z_j} =:\theta_j, $$
we can write the vector field in the form
\[ X= \sum_{j=1}^{n} b_j(z) z_j\d_{z_j}, \;\;\;\;b_j(z) \in A=\CM\{z,z^{-1}\}.\]
In particular, the constant vector field
\[\sum_{j=1}^n \nu_j \d_{x_j}\]
becomes, up to a factor $2 \pi i$, the linear vector field 
\[ \sum_{j=1}^n \nu_j z_j \d_{z_j}=\sum_{j=1}^n \nu_j \theta_j\]
in the $z$-variables.
\subsubsection{}
We denote the set of all analytic vector fields by $\Theta(A)$ and we
have an isomorphism of $A$-modules 
\[ A^n \stackrel{\simeq}{\to} \Theta(A),\;\;(a_1(z),\ldots,a_n(z)) \mapsto \sum_{j=1}^n a_j(z) \theta_j\]
and algebraically one may identify $\Theta(A)$ with the module of {\em derivations} of the ring $A$:
\[ \Theta(A) \stackrel{\approx}{\to}\Der_{\CM}(A),\;\;\; X \mapsto (f \mapsto X(f)).\]
The commutation of derivations gives $\Theta(A)$ a
natural structure of a Lie-algebra.

\subsubsection{}
As the torus $\TM^n$ is compact, a real analytic vector field $X$ 
has a globally defined flow 
\[\Phi_t: \TM^n \to \TM^n,\;\;\;t\in \RM,\]
consisting of real analytic automorphisms of the torus. These automorphisms
act on the ring $A=\CM\{z,z^{-1}\}$ by composition.
This action can be described formally as the exponentiation of the vector field:
\[ f \mapsto f\circ \Phi_t=e^{tX}(f)=f+tX(f)+\frac{t^2}{2!}X(X(f))+\ldots\]
In section \ref{S::Kspaces} we  give a functional analytic treatment of the exponential
of a vector field, which will provide us with explicit norm-estimates, which
leads to a direct proof of the convergence of the series and the existence 
of the flow.
\subsubsection{}
The above exponential series $e^{tX}$ defines an automorphism of $A=\CM\{z,z^{-1}\}$. There is a corresponding induced adjoint action on the
 module $\Theta(A)$ of derivations
\[ Y \mapsto \Phi_t \circ Y \circ \Phi_t^{-1},\] 
that we denote by
\[ e^{[tX,-]}=1+t[X,-]+\frac{t^2}{2!}[X,[X,-]]+\ldots \in Aut(\Theta(A)),\]
so that
\[ e^{[tX,-]}(Y)=Y+t[X,Y]+\frac{t^2}{2!}[X,[X,Y]]+\ldots\]
\newpage
\subsection{Symplectic vector fields on the torus}

\subsubsection{}
Via the embedding $\TM:=\TM^{2d} \subset (\CM^*)^{2d}$ the symplectic form 
\[\omega=\sum_{1 \le i<j \le 2d} \omega_{ij}dx_i \wedge dx_j\] 
on $\TM$ can be seen as the restriction of the (complex) holomorphic 
symplectic form
\[\frac{1}{(2 \pi i)^2}\sum_{1 \le i <j \le 2d} \omega_{ij} \frac{dz_i}{z_i}\wedge \frac{dz_j}{z_j}\]
on $(\CM^*)^{2d}$. Such a choice of symplectic form gives the ring
$A=\CM\{z,z^{-1}\}$ the structure of a {\em Poisson-algebra}, with Poisson-bracket
\[\{f,g\}:=\sum_{i,j=1}^{2d} \omega_{ij} (\theta_i f) (\theta_j g)\]
which coincides, up to a factor $(2\pi i)^2$, with the usual Poisson-bracket
\[    \sum_{i,j=1}^{2d} \omega_{ij}(\d_{x_i} f)(\d_{x_j}g) , \]
when written in the original variables $x_1,x_2,\ldots,x_{2d}$.

\subsubsection{}
The (complex) symplectic vector fields on $\TM$ are denoted by
\[ \St(A) \subset \Theta(A)\]
and are in one-to-one correspondence with closed one-forms. In particular the 
interior product of the symplectic form $\omega$ with the fields 
$\theta_i=z_i \d_{z_i}$ give constant and therefore 
non-exact closed one-forms 
\[ \a_i:=\iota_{\theta_i}(\omega),\;\;\;i=1,2,\ldots,2d .\]
As
$$\d_{x_1} \w \d_{x_2} \w \dots \w \d_{x_{2d}} \neq 0, $$
these one-forms generate the De Rham cohomology group $H^1_{dR}(\TM,\CM)$
of the torus. Consequently, any closed $1$-form $\a$ can be written as
$$\a=\sum_{i=1}^{2d} c_i\a_i+dh . $$
From the dual perspective this means that any symplectic vector field $S$ 
can be written as the sum of a quasi-periodic field and a Hamiltonian part:
$$S=\sum_{i=1}^{2d} c_i \theta_i+\{-,h \}. $$
Note that this representation is essentially unique; the function $h$ is
determined up to a constant. We will always choose $h$ to have {\em vanishing
constant term} when written as an analytic Fourier series.

\subsubsection{}
One can summarise the above discussion by saying that the cokernel
of the map 
\[ A  \to \St(A),\;\;\; f \mapsto \{-,f\}\]
is identified with the homology group $H_1(\TM,\CM)$, whereas the kernel
consists of the constants. Hence there is a natural exact sequence
\[ 0 \to \CM \to A \to \St(A) \to H_1(\TM,\CM) \to 0\]
and the image of the map in the middle consists precisely of the 
Hamiltonian vector fields.
 
%%%%%%%%%%%%%%%%%%%%%%%%%%%%%%%%%%%%%%%%%%%%%%%%%%%%%%%%%%%%%%%%%%%%%%%%%%%%%%
\subsection{The torus with formal parameters}
We will pick a reference symplectic form $\omega^0$ and add parameters 
that describe the variation of the symplectic form on the torus. Furthermore
we add parameters that detune the frequencies. The reason for this is that it
makes it easier to formulate the normal form iteration. By the implicit function
theorem we will later express the perturbed frequencies in terms of the 
perturbed symplectic form. However, initially we will consider these as formal,
independent parameters.

\subsubsection{ The ring $R$}

We add variables 
$$\phi_1,\dots,\phi_{2d}$$
to parametrise the frequencies and
$$\dt_{1},\dt_{2},\dots,\dt_{l},\ l=d(2d-1) ,$$ 
to parametrise symplectic forms in the neighbourhood of $\omega^0 \in \SM$ 
and set
$$\omega_\dt:=\omega^0+\sum_{i}\dt_{k}\omega_k $$
where $\omega_k:=dx_{i_k} \w dx_{j_k}.$
The elements of the ring
$$R:=\CM\{z,z^{-1} \}[[\dt,\phi]] $$
are formal series
\[ f=\sum_{K,L} f_{K,L} \dt^K\phi^L ,\]
where the coefficients $f_{K,L} \in \CM\{z,z^{-1}\}$.
We assign the weight one to the variables $\dt_i$ and $z_i$ 
and weight zero to $\phi_i$. We can naturally extend the
truncation from the ring $A$ to the ring $R$, so that
\[ [f]^k :=\sum_{L}\sum_{|K| \le k} [f_{K,L}]^{k-|K|} \dt^K \phi^L.\]
\subsubsection{}
The ring $R$ has a natural Poisson structure that can be written as
$$\{ f,g \}:=\sum_{i < j} (\omega_{\dt})_{ij}(\theta_i f)\wedge (\theta_j g) . $$
The sub-ring 
$$R_0:=\CM[[\dt,\phi]] $$ 
is the centre of this Poisson algebra and Poisson vector fields are of the form:
$$v+ \sum_{i=1}^{2d} c_i \theta_i+\{-,h\},\; v \in \Der_\CM(R_0),\;\;c_i \in R_0,\;\; h \in R .$$
Here $ \Der_\CM(R_0)$ is the free $R_0$-module spanned by the $\partial_{\phi_i}$ and the $\partial_{\dt_i}$. For a symplectic vector fields one has $v=0$.

We have inclusions of Lie algebras:
$$\St(R) \subset \Pt(R) \subset \Der_\CM(R) $$
where $\St$ and $\Pt$ stand respectively for symplectic and Poisson vector fields.
There is an intermediate space of Poisson vector fields $\Qt(R)$, where the
vector field $v$ only involves derivations with respect to $\phi_i$:
$$\Qt(R):=\St(R) \bigoplus_{i=1}^{2d} R_0\d_{\phi_i}  $$
Clearly,
\[ \St(R) \subset \Qt(R) \subset \Pt(R)\]
and $\Qt(R)$ is a Lie-subalgebra of $\Pt(R)$. This sub-Lie-algebra $\Qt(R)$
will be relevant in our iteration process.

\section{The iteration in the ring $R$}
In this section we outline the iteration process that brings a perturbed
vector field back to normal on the level of power series.
In later sections we lift this iteration to the level of Banach spaces and 
give the estimates that lead to a proof of convergence. We start with a discussion of the associated {\em homological equations}, i.e. the linear operators 
that we will need to invert.
\newpage
\subsection{The homological equation for functions}
\subsubsection{}
One particular vector field that will play a crucial role in the sequel
is the following. For a quasi-periodic vector field
\[X_{\nu}:=\sum_{i=1}^{2d} \nu_i(\dt) \theta_i\]
with non-resonant frequency vector 
\[\nu(0)=(\nu_1(0),\nu_2(0),\ldots,\nu_{2d}(0))\]
we consider the {\em versal unfolding} 
\[ \Vt:=X_{\nu}+\sum_i \phi_i \theta_i=\sum_{i=1}^{2d}(\nu_i(\dt)+\phi_i) \theta_i\]
of $X_{\nu}$ obtained by {\em detuning} the frequencies with the variables 
$\phi_1,\phi_2,\ldots,\phi_{2d}$ of the ring $R$.
We study the action of this vector field on functions and vector fields.
\subsubsection{}
The Lie-derivative
$$D:R \to R,\ f \mapsto L_{\Vt} f =\Vt(f)$$
is diagonal in the monomial basis:
$$D: z^I \mapsto (\nu(\dt)+\phi,I) z^I , $$
where $(-,-)$ denotes the Euclidean scalar product.\\
Hence the operator $D$ is equal to taking Hadamard product with
$$g(z):= \sum_{I \in \ZM^{2d}}(\nu(\dt)+\phi,I) z^I.$$
so that
$$ D(f) = g \star f .$$
The kernel of the map $D$ is the centre $R_0=\CM[[\dt,\phi]]$ of the Poisson algebra $R$.
\subsubsection{}
As we assumed $\nu(0)$ to be non-resonant, the functions 
\[(\nu(\dt)+\phi,I),\;\; I \neq 0,\] are invertible elements of 
$R_0=\CM[[\phi,\delta]]$, and thus we can consider the formal power series:
$$H(z):= \sum_{I \in \ZM^{2d} \setminus \{ 0 \}}(\nu(\dt)+\phi,I)^{-1} z^I \in R_0[[z,z^{-1}]] ,$$
which we call the {\em resolvente} of $\Vt$. Note that if we try to interprete
this formal series as a function, it would have, in general, poles of order 1
along a dense set of hyperplanes defined by the resonance conditions
$(\nu(\dt)+\phi,I)=0$, $I \in \ZM^{2d}$.

If $P$ is a Laurent polynomial in $z$ {\em with vanishing constant coefficient}, then
$$D( H \star P)= g \star H \star P=P , $$
so we get an inverse to the operator $D$.
If we want to write a similar formula for more general elements of $R$, we
need to assume that the frequency vector $\nu(0)$ satisfies 
a {\em Diophantine condition $C,N$}:
\[ |(\nu(0),I)| \ge \frac{C}{|I|^N}.\]
It is readily shown that under such a Diophantine condition on $\nu(0)$
the operator $H \star$ maps $R$ to itself, and provides an inverse 
to the Lie-derivative operator.
$$ D( H \star f)= f,\;\;\; f\in R $$
We note that the Lie-derivative $L_{\Vt}$ and therefore the map $D$ extends 
to arbitrary tensors. The Lie-derivative commutes with the exterior 
derivative and therefore
$$Dd  \left( H \star P \right)=dD \left( H \star P \right)=d  P .$$
%%%%%%%%%%%%%%%%%%%%%%%%%%%%%%%%%%%%%%%%%%%%%%%%%%%%%%%%%%%%%%%%%%%%%%%%%%%%%%%%
\subsection{The homological equation for vector fields}

\subsubsection{}
The action of the Lie-derivative on vector fields is given by the
Lie bracket or commutator of vector fields:
\[ L_{\Vt} Y = [\Vt,Y]=-L_Y \Vt .\]

For a given symplectic vector field $S$, we want to solve the 
{\em homological equation}
$$L_{\Vt} Y  =  S $$
for the vector field $Y$. Somewhat surprisingly, this can be done 
quite explicitly.

\subsubsection{} If we assume a Diophantine condition on $\nu(0)$
we have the following: 

\begin{proposition}\label{homolog}
Decompose $S$ into Hamiltonian and non-Hamiltonian part:
$$ S=\{-,f\}+\sum_{i=1}^{2d} c_i \theta_i,\;\;\;c_i \in R_0$$
Then the equation
$$L_{\Vt} Y  =  S $$
is solved by
$$Y:=\{-,H \star f \}+\sum_{i=1}^{2d} c_i \d_{\phi_i} , $$
where $H=H(z)$ is the resolvente introduced above.
\end{proposition}

\begin{proof}
Write $D=L_{\Vt}$ for the operation of Lie-derivative with respect to $\Vt$.
As the symplectic form $\omega_{\dt}$ is invariant under the flow of $\Vt$, one
has $D \omega_{\dt}=0$ and consequently
\[ D(\{f,g\})=\{Df,g\}+\{f,Dg\}.\]
The commutator of the vector fields $\Vt$ and $\{-,g\}$ acting on $f$
is 
\[[\Vt,\{-,g\}](f)=\Vt(\{f,g\})-\{\Vt(f),g\}=D(\{f,g\})-\{Df,g\}=\{f,Dg\},\]
which means that
\[ D(\{-,g\})=\{-,Dg\}.\]
So we get
\begin{align*}
DY&=\{-,D(H \star f) \}+[\Vt,\sum_{i=1}^{2d} c_i \d_{\phi_i}] \\ 
 &=\{-,f \}+\sum_{i=1}^{2d} c_i \theta_i \\
 &=S
\end{align*}
\end{proof}

\subsubsection{}
The decomposition of $S$ into Hamiltonian and non-Hamiltonian
parts is unique. Therefore, still under a Diophantine 
condition on the frequency $\nu(0)$,  the explicit solution to 
homological equation determines a map
\[
\begin{array}{r c l} j: \St(R) &\to & \Qt(R)\\
  S &\mapsto & Y=\{-,H \star f \}+\sum_{i=1}^{2d} c_i \d_{\phi_i} .  
\end{array}
\]

\subsubsection{}
If $\nu(0)$ is not Diophantine, then taking the Hadamard product with the
resolvente $H$ no longer maps $R$ to itself. As a result, the
map $j$ defined above does not make sense. However, it turns out
that it is sufficient to consider an {\em approximate inverse}.
If $S=\{-,f\}+\sum_{i=1}^{2d} c_i \theta_i$ and $A\subset \ZM^n$
a finite subset, we set
\[[S]_A:=\{-,[f]_A\}+\sum_{i=1}^{2d} [c_i]_A \theta_i \]
and can define maps
\[ j_A :\St(R) \to \Qt(R)\]
by setting
\[j_A(S):=j([S]_A)=\{-, H \star [f]_A\}+\sum_{i=1}^{2d} [c_i]_A \d_{\phi_i} .\]

We also use the notations
\[ [S]^m:=\{-,[f]^m\}+\sum_{i=1}^{2d}[c_i]^m\theta_i,\;\;[S]_m:=\{-,[f]_m \}+\sum_{i=1}^{2d}[c_i]_m\theta_i  \]
for the truncations of a vector field in lower and higher Fourier modes.

%%%%%%%%%%%%%%%%%%%%%%%%%%%%%%%%%%%%%%%%%%%%%%%%%%%%%%%%%%%%%%%%%%%%5%%%%
\subsection{The iteration}

\subsubsection{}
We start with a quasi-periodic vector field
$$X_{\nu}= \sum_{i=1}^{2d} \nu_i(\dt) z_i \d_{z_i}$$
and a perturbation $X$ of it by a symplectic vector field $S=S_0 \in \St(R)$:
\[X:=X_{\nu}+S_0 .\]
Furthermore, the corresponding unfolding of $X_{\nu}$ 
$$\mathcal{V}:=X_{\nu}+\sum_{i=1}^{2d} \phi_i \theta_i=\sum_{i=1}^{2d} 
(\nu_i(\dt)+\phi_i)\theta_i$$ 
is obtained by introducing the detuning parameters $\phi_i$. In a similar way,  
the perturbed field $X$ can be detuned to
\[X_0=X+\sum_{i=1}^{2d}\phi_i \theta_i=\Vt+S_0 .\]

\subsubsection{}
We decompose the perturbation into two parts:
\[ S_0=[S_0]^{2}+[S_0]_{2}.\] 
Using proposition \ref{homolog} we can find $Y_0 \in \Qt(R)$ that solves
the equation:
\[
[Y_0,\mathcal{V}]=[S_0]^{2}
\]
The vector field $Y_0$ integrates to an automorphism $e^{-Y_0}$ of $R$.
By the adjoint action, the vector field $X_0$ is transformed into
$$X_1:=e^{-[Y_0,-]} X_0=X_0-[Y_0,X_0]+\ldots$$
and define the next perturbation $S_1$ by setting
$$ X_1= \mathcal{V}+S_1 .$$

\subsubsection{}
We repeat this operation by taking terms up to order $2$, $2^2$, $2^3$ and 
so on. We obtain an iteration scheme of the form
\begin{align*}
X_n&=\mathcal{V}+S_n, \\
S_n&=[S_n]^{2^{n+1}}+[S_n]_{2^{n+1}},\\
[Y_n,\mathcal{V}]&=[S_n]^{2^{n+1}},\\
X_{n+1}&=e^{-[Y_n,-]}X_n.
\end{align*}

\subsubsection{}
The iteration may at first look foolish from a formal point of view: 
as the vector fields $Y_n$ will contain all terms of degrees up to $2^n$, its
exponential will reintroduce monomials that one tries to remove. 
One is reminded of the mythos of {\sc Sisyphos}, but we will show later by 
a direct estimate that the norm of the remainder decreases quadratically 
in appropriate Banach spaces, so that his burden is descreasing quickly, 
although the removal of even the first term will require an infinite numer 
of iterations and keeps him busy forever.

%%%%%%%%%%%%%%%%%%%%%%%%%%%%%%%%%%%%%%%%%%%%%%%%%%%%%%%%%%%%%%%%%%%%%%%%%%%%%%
\subsection{Relation to the original vector field}
We start from a perturbed $X=X_{\nu}+S_0$ quasi-periodic motion 
and considered the corresponding perturbation of the versal unfolding 
\[X_0=\Vt+S_0.\] 
The iteration produces a sequence of vector fields
\[ Y_0,\;Y_1,\;Y_2,\;\ldots,Y_n,\;\ldots\] 
and by exponentiation we obtain corresponding automorphisms of $R$
$$\p_n:=e^{-Y_n}\dots e^{-Y_0}  \in Aut(R)$$
By the corresponding adjoint automorphism 
\[\psi_n:=e^{[Y_n,-]}\ldots e^{[Y_0,-]} \in Aut(\Theta(R))\]
we obtain a corresponding sequence of transformed vector fields 
\[ X_0,X_1,\dots,X_n,\dots\] 
with
$$X_n=\mathcal{V}+S_n=\psi_n(\mathcal{V}+S_0) .$$

Trivially, our original vector field $X$ is obtained from $X_0$ by setting to zero all the detuning variables $\phi_k=0,\;\;k=1,2,\ldots,2d$:
\[ X=(X_0)_{\phi=0}.\]
The automorphisms $\p_n \in Aut(R)$ clearly have to preserve the Poisson-center, so induce automorphisms
\[\p_n: R_0 \to R_0 .\]
By construction, the vector fields $Y_n$ only contain $\partial_{\phi_i}$, and
not the $\partial_{\dt_i}$, so the automorphisms $\p_n$ have the additional
property that 
\[ \p_n(\dt_i)=\dt_i,\;\;i=1,2,\ldots,l\,\]
and hence map the element $\phi_k$ to certain power series
\[ \p_n(\phi_k) =:R_{n,k}(\phi_1,\ldots,\phi_{2d},\dt_1,\ldots,\dt_l)=:R_{n,k}(\phi,\dt),\;\;\;k=1,2,\ldots,2n.\]
As 
\[R_{n,k}(\phi,\delta)=\phi_k+O(2),\]
we may apply the (formal) implicit function theorem and solve the $\phi_k$
uniquely in terms of the $\dt$'s.
\[ \phi_k=g_{n,k}(\dt).\]
The functions $g_{n,k}$ describe how we have to adapt the frequency to the
parameters that change the symplectic form. Under the adjoint automorphism 
$\psi_n$, the vector field $X_0+S_0$ is mapped to $X_n+S_n$. As the condition
$\phi_k=0, k=1,2,\ldots,2d$ is transformed into $R_{n,k}(\phi,\dt)=0,k=1,2,\ldots,2d$, we see:\\

\begin{proposition} The adjoint automorphism $\psi_n$ transforms the vector 
field 
$X_{\nu}+S_0$ into
\[\psi_{n}(X_{\nu}+S_0)=(\Vt)_{\phi_k=g_{n.k}(\dt)}+(S_n)_{\phi_k=g_{n,k}(\dt)}.\]
\end{proposition}

So in the limit $n \to \infty$, we expect to get a relation of the form
\[\psi_{\infty}(X_{\nu}+S_0)=(\Vt)_{\phi_k=g_{\infty,k}(\dt)}=\sum_{i=1}^{2d}(\nu_i(\dt)+g_{\infty,i}(\dt))\theta_i,\]
hence we produce a coordinate transformation that conjugates the
symplectic perturbation $X_{\nu}+S$ of a quasi-periodic vector field $X_\nu$
to a nearby quasi-periodic vector field, with frequency that depends on the
perturbation of the symplectic form. 
%%%%%%%%%%%%%%%%%%%%%%%%%%%%%%%%%%%%%%%%%%%%%%%%%%%%%%%%%%%%%%%%%%%%%%%%%%%%

\subsection{Almost quadratic nature of the iteration}
\subsubsection{}
Our iteration is defined by first writing
\[ X_n=\Vt+S_n\]
and then recursively
$$\left\{
\begin{array}{rcl}
Y_n&=&j_n(S_n) \\
X_{n+1}&=&e^{-[Y_n,-]}X_n
\end{array}
\right.
$$

Here we use the notation 
\[ j_n(-)=j_{A_n}(-)\]
where $A_n$ is the set of monomials, whose absolute value of weight
is smaller than $2^{n+1}$.

\subsubsection{}
The iteration is quadratic with remainder term in the following sense:
\begin{align*}
X_{n+1}&=e^{-[Y_n,-]}\left(\Vt+[S_n]^{2^n} \right)+e^{-[Y_n,-]}\left([S_n]_{2^n} \right) \\
&=e^{-[Y_n,-]}\left( \Vt+[Y_n, \Vt]\right)+e^{-[Y_n,-]}\left([X_n]_{2^n} \right)\\
&=\Vt+(e^{-[Y_n,-]}(\Id+[Y_n,-])-\Id)\Vt+e^{-[Y_n,-]}\left([X_n]_{2^n} \right)
\end{align*}
For a power series in a single variable $x$ 
\[f(x)=\sum_{i=1}^\infty a_i x^i\]
and a vector field $X$, we put
\[f_*(X):=\sum_{i=1}^{\infty}a_i (L_X)^i\]
So if we write
$$f(x)=e^{-x}(1+x)-1=-x^2+o(x^2) ,$$
the iteration can be written in the form
$$S_{n+1}=f_*(j_n(S_n))\Vt+ e^{-[j_n(S_n),-]}\left([S_n]_{2^n} \right)$$

So in a formal sense the iteration has a quadratic term 
$f_*(j_n(S_n))$ and a remainder part  $e^{-[j_n(S_n),-]}([S_n]_{2^n})$.
Although terms of low degree remain at each step of the iteration, it might 
be expected that, because of this quadraticity, their coefficients rapidly 
tend to zero. We will see that this is, under certain conditions, indeed the 
case.

\subsubsection{}
The above can be seen as an iteration in the Fr\'echet space $\St(R)$.
We will now formulate a version of the iteration in terms of a system
of Banach spaces of holomorphic functions attached to neighbourhoods $T_r$ 
of $\TM$ in $(\CM^*)^{2d}$. We will keep track of the norms of $S_n$ 
during the iteration. This will show that the above process converges 
over an non-trivial Cantor-like set, defined by the condition that the 
norm remains sufficiently small.

%%%%%%%%%%%%%%%%%%%%%%%%%%%%%%%%%%%%%%%%%%%%%%%%%%%%%%%%%%%%%%%%%%%%%%%%%%%%%
\section{A short review on functorial analysis}
\label{S::Kspaces}
In the proof we have to work with many different Banach spaces and various maps
between them that are 'compatible' in various ways. These maps result from 
restriction maps that appear in shrinking of domains during the iteration 
process, or changes in the type of Banach space considered. 
Of course, one needs to have explicit control over all the norms of these maps. 
To keep track of all these, we found it convenient to use an abstract framework 
that was developed in \cite{KAM_book} to which we refer for more details.
 
\subsection{Kolmogorov spaces}
We give a quick overview of the formalism of Banach spaces parametrised by
ordered sets that one encounteres often in dealing with function spaces over
shrinking domains of definition.

\subsubsection{}
We denote by {\bf Ban} the category whose objects are Banach spaces, and whose
morphisms are bounded linear operators. If $B$ is a small category, we
mean by a {\em Banach space over B} a covariant functor
$$F:B \to {\bf Ban} .$$ 
The {\em total space} of such a Banach space over $B$\footnote{We denote by $B$ also the set objects of the categroy $B$} is defined as
$$E:=\bigsqcup_{b \in B} E_b,\ E_b=F(b)  $$
and there is a natural map
$$E \to B $$
which maps the elements of $E_b$ to $b$. We sometimes use the notation $(b,x), x \in E_b$ for the elements of $E$ and we use the generic name $|-|_b$
for the norm on the Banach space $E_b$. For an element $x_b \in E_b$ we often 
write $|x_b|$ instead of $|x_b|_b$, etc. A {\em section over $A \subset B $} is
a choice of vectors $x_b \in E_b$ for all $b \in A$, like in the theory of vector bundles or sheaves.
A section over $A$ is called {\em bounded} if the function $b \mapsto |x_b|$ is bounded on $A$.\\

\subsubsection{} 
The functor property means that for each morphism from $a$ to $b$ in $B$, there is a corresponding continuous linear map $E_a \to E_b$ that is part of the structure. But we will often just say that $E$ is a Banach space over $B$, the rest of the structure being understood.

Any {\em partially ordered set}
%\footnote{We use the french convention {\em ordered} means that there is a {\em partial order} .}
%$(B,\ge)$ 
is naturally a small category with spaces of morphism $Mor(t,s)$ consisting 
of a single element if $t \ge s$. 
So a Banach space over $(B,\ge)$ consist of Banach spaces $E_t, t \in B$ and
for each $t \ge s$ compatible continuous linear maps
\[ e_{st}: E_t \to E_s,\]
called {\em restriction maps}, where compatibility means
\[e_{us}\circ e_{st}=e_{ut}\;\; \textup{if}\;\;t \ge  s \ge  u .\]
All examples that we will consider in this paper are of this type. 

Using these maps we can compare the spaces over different points
and a section $x_t$ over a subset $A \subset B$ 
is called {\em horizontal} if it is {\em compatible} with the 
restriction maps in the sense that
\[e_{st}(x_t)=x_s , \;\; t \ge s,\;\; t,s \in A.\]

\subsubsection{}
A Banach space over $B$ is called a {\em Kolmogorov space}, if all the 
restriction maps
$$ e_{st} \in Hom(E_t,E_s)  $$
have norm $\le 1$, where we put the operator norm on the space $Hom(E_t,E_s)$ of
continuous linear maps. The Banach spaces 
\[ E_t:=C^0([0,t],\RM)\] 
of continuous functions on the interval $[0,t]$, $t \in ]0,\infty[$,
with the supremum norm and obvious restriction mappings 
\[E_t \to E_s,\;\; t \ge s,\]
is a good example of a Kolmogorov space over $]0,\infty[$ to keep in mind.\\

Kolmogorov spaces over an interval $]0,S]$ occur frequently and we refer to
them as $S$-Kolmogorov spaces for short. More generally a Kolmogorov space 
is called a $Kn$-space if the base $B$ is a subset 
of $(\RM^n, \ge )$ where the partial order is taken {\em component wise}:
$$x \ge y \iff  x_i \ge y_i\;\; \textup{for}\;\;i=1,2,\ldots,n .$$ 
In this paper we will encounter only $K1$ and $K2$ spaces.

\subsubsection{}
Let $B$ be a partially ordered set, considered as a category.
To any Banach space $E$ over $B$ we can associate in a natural
way a Kolmogorov space $EK$ over $B$, by setting
\[ EK_b :=\Gamma^{\infty}(]-\infty,b],E) \]
which is the space of {\em bounded horizontal sections} of $E$ over the
{\em down-set} of $b$: 
\[ ]-\infty,b] :=\{ b' \in B\;|\; b' \leq b\}\]

If $x=(x_t),\;\;t \le b$, is such a section, we assign to it the norm
\[ |x|_b:=\sup_{t \le b } |x_t|\]  
Clearly, if $ b' \le b$, then we have $]-\infty,b']\, \subset\,\, ]-\infty,b]$, so that
indeed
$$|x|_{b'} \leq |x|_{b},$$
as we are taking the supremum over a smaller set and so
the natural restriction mappings $EK_b \to EK_{b'}$ have norm $ \leq 1$.

\begin{proposition}[\cite{KAM_book}]
Given a relative Banach space $E$ over $B$, the associated space 
$EK$ over $B$ is a Kolmogorov space.
\end{proposition}
\begin{proof}
The only non-trivial fact to check is the completeness of the normed vector 
spaces $EK_b$. So let $\g_n$ be a Cauchy sequence of sections in $EK_b$. 
 
For any $b'>b$, the sequence $\g_n(b')$ is a Cauchy sequence in $E_{b'}$ and
therefore converges to a limit $\g(b')$. We need to show that the norm of this limit is finite.

The norms $|\g_n| $ form a Cauchy sequence of real numbers  and therefore converges to a limit $M$.
\begin{align*}
|\g(b')| &\leq |\g(b')-\g_n(b')|+|\g_n(b')| \\
          & \leq |\g(b')-\g_n(b')|+M .
\end{align*}
So, passing to the limit, we see that the norm of $\g$ is bounded by $M$ and in fact equal
to $M$ (because the norm is a continuous map). 
\end{proof}
For this reason we call $EK \to B$ the {\em Kolmogorification} of
$E \to B$; if $E$ is already Kolmogorov, then
$EK=E$, so we get back the original space\footnote{The process of Kolmogorification is somewhat analoguous to sheafification of a presheaf.}.

\subsubsection{}
Given a Kolmogorov space 
$$E \to B $$
we can form in a natural way an {\em opposite Kolmogorov space}
$$E^{op} \to B^{op} $$

The underlying set $B^{op}$ is the same as $B$, but the order relation
on $B^{op}$ is {\em opposite} to that of $B$: it $t \ge s$ in $B$, then $s \ge t$ in $B^{op}$. 

The fibre at $b$ is
$ E^{op}_b:=\G^\infty([b,+\infty[,E),$ 
the Banach space of horizontal bounded sections over the {\em up-set} of
$b$
$$[b,+\infty[=\{b' \in B\;|\; b'\ge b\}$$ of $b$, with the supremum norm. 
So the opposite space of a Kolmogorov space is again a Kolmogorov space.

%%%%%%%%%%%%%%%%%%%%%%%%%%%%%%%%%%%%%%%%%%%%%%%%%%%%%%%%%%%%%%%%%%%%%%%%%%%%%%
\subsection{Local operators}

\subsubsection{}
We will consider the case where $B:=]0,S]$ is an interval in $\RM$, 
with the natural ordering $\ge$, so we are dealing with $S$-Kolmogorov
spaces. If $E$ and $F$  are $S$-Kolmogorov spaces, then spaces 
$\Hom(E_t,F_s)$ form the fibres of a K2-space 
\[{\Ht}om(E,F) \to B^{op} \times B\] 
in a natural way. The restriction maps from 
$\Hom(E_t,E_s)$ to $\Hom(E_{t'},F_{s'})$ are obtained 
as composition with the restriction maps of $E$ and $F$:
\[ u_{s't'}=f_{s' s}u_{st}e_{tt'},\]
which forces to have $ s \ge s'$, but $t' \ge t$. By restriction to the
triangle
$$\D:=\{(t,s) \in B^2:t>s \} \subset B^{op}\times B$$ 
and a rescaling of the norm we obtain the K2-Kolmogorov space 
$${\Ht}om^k(E,F) \to \D,$$
whose fibre is the space $\Hom(E_t,F_s)$ of continuous linear mappings
with {\em rescaled norm}
$$|t-s|^k\,\| u_{st} \| .$$
Here and in the sequel $\| \cdot \|$ stands for the operator norm:
\[ \| u_{st} \|:=\sup_{x \in E_t \setminus \{ 0 \}}\frac{| u(x)|_s}{| x |_t}.\]

\subsubsection{}
We consider the opposite space
$${\Ht}om^k(E,F)^{op} \to \D^{op} ,\;\; \D^{op} \subset B \times B^{op}$$
This relative Banach space defines, after Kolmogorification, a Kolmogorov 
space denoted by $\Lt^k(E,F)$ and called the {\em $K2$-space of $k$-local operators.}
Unraveling the definition, we have
\[\Lt^k(E,F)_{s,t} := \Gamma^{\infty}(\Delta(t,s),{\Ht}om^k(E,F)),\]
whose elements are {\em compatible systems} $u=(u_{a,b})$
\[ u_{a,b} \in \Hom(E_b,F_a),\;\;\; (b,a) \in \Delta(t,s),\]
where
$$\D(t,s):=\{(b,a) \in B^2:t \geq b>a \geq s  \} \subset B^{op}\times B$$ 
for which 
\[|u|_{s,t}:=\sup_{(b,a) \in \D(t,s)} |b-a|^k\,\| u_{a,b} \| <\infty. \]

\begin{center}
\includegraphics[height=5cm]{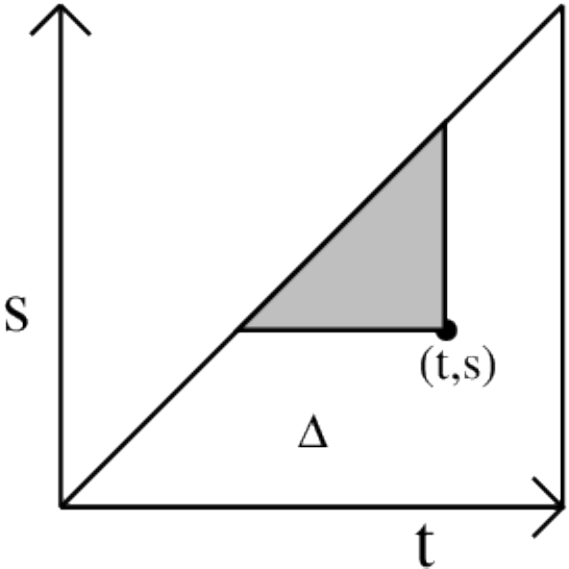}
\end{center}

As explained in \cite{KAM_book}, one can form the {\em direct image} of $\Lt^k(E,F)$ under the projection
$$\pi:(t,s) \mapsto t $$
and obtain a K1 space:
 $$L^k(E,F) \to B, $$
which we call it the {\em $K1$-space of $k$-local operators}. 
An element of the Banach space $L^k(E,F)_t$ is simply
a {compatible} family $(u_{r,s})$ of Banach space morphisms, defined for 
all  $r < s \le t$, such that
$$| u|_{t}:=\sup_{r<s \leq t} |s-r|^k\,\| u_{r,s} \|  $$
is finite. 

Our formalism provides an automatic method of handling such compatible 
families of operators, and packing them into Banach spaces with appropriate 
norms in a systematic way.

\subsubsection{}
If $E,F,G$ are Kolmogorov spaces over $B$, we have a well 
defined bilinear composition map
\begin{align*}
\circ: \Hom_{\D(t,s)} (E,F)) \times \Hom_{\D(t,s)} (F,G)) &\to \Hom_{\D(t,s)}(E,G))\\
 (u,v)& \mapsto v \circ u \end{align*}
which we leave to the reader to spell out (see also \cite{KAM_book}).
Here we use to notation
\[ \Hom_A(E,F):=\Gamma^{h}(A,{\Ht}om(E,F))\]
for the space of horizontal sections of ${\Ht}om(E,F)$
over $A$. The sub-spaces of $k$-local operators
\[ \Lt^{k}(E,F)_{t,s} \subset \Hom_{\D(t,s)} (E,F))\]
behave well under this composition:

\begin{proposition}
Let $E,F,G$ be $S$-Kolmogorov spaces.
If $u \in \Lt^{k}(E,F)$ and $v \in \Lt^{m}(F,G)$ then 
$$v\circ u \in \Lt^{k+m}(E,G) .$$ 
 Moreover, one has the norm estimate: 
\[ |v \circ u| \le  \frac{(k+m)^{k+m}}{k^{k} m^{m}} |v||u|.\]
\end{proposition}
\begin{proof} For  $s < s' < t $ we have
\[ \|(v \circ u)_{st}\| \le \|v_{ss'}\|\,\|u_{s't}\| .\]
 From the definition of locality we have:
\begin{align*}
\|v_{ss'}\| &\le \frac{|v|}{(s'-s)^{k}},\\
\|u_{s't}\| &\le \frac{|u|}{(t-s')^{m}}
.\end{align*}
We take the point $s'$ such that 
$$(s'-s)=\frac{m}{k+m}(t-s),\ (t-s')=\frac{k}{k+m}(t-s) ,$$
and  find the estimate
$$ \|v_{ss'}\|\,\|u_{s't}\| \leq \frac{(k+m)^{k+m}}{k^{k} m^{m}} \frac{|v| |u|}{(t-s)^{k+m)}} ,$$
thus
\[ |v \circ u| \le \frac{(k+m)^{k+m}}{k^{k} m^{m}} |v| |u| .\]
\end{proof}

The form of the above estimate leads to the idea to {\em recalibrate} the norm on $\Lt^{k}(E,F)$
and introduce
\[ ||u||:=\left(\frac{e}{k}\right)^k |u|,\;\;\;e=2.718281828\ldots\]
This recalibrated norm now is {\em sub-multiplicative}:
\[ ||v \circ u|| \le ||v||\, ||u|| .\]
For this to work one could replace Eulers constant $e$ in the above definition by any positive number, but this choice has some advantages.  

\begin{corollary}\label{C::powerestimate} If $u_1,u_2,\ldots,u_n \in \Lt^1(E,E)$ then
\[ u_1 \circ u_2\circ \ldots \circ u_n \in \Lt^n(E,E),\] 
and we have the estimate:
\[ ||u_1 \circ u_2\circ \ldots \circ u_n|| \le ||u_1|| ||u_2||\ldots||u_n||, \]
or, in terms of the original norm:
\[ |u_1 \circ u_2 \circ \ldots \circ u_n| \le n^n |u_1| |u_2| \ldots |u_n| .\]
Especially for $u_1=u_2=\ldots=u_n=u$:
\[ ||u^n|| \le ||u||^n,\;\;\; |u^n| \le n^n |u|^n .\]
\end{corollary}

In reading these formulas one should be aware that for reasons of simplicity
we always write $|-|$ or $||-||$ for the norms, but which norm this in fact is, 
depends on the space in which the element we apply the norm to is considered. 
For example,  if $u$ is 1-local, the  norm in the expression $|u|^n$ refers 
to the $n$-th power of the norm on $\Lt^1(E,E)$, whereas in $|u^n|$ we are 
using the norm on $\Lt^n(E,E)$~! In practice, this notational ambiguity always 
is resolved in a unique way.

\subsubsection{}
There is a canonical $0$-local map $\iota$ of a Kolmogorov space to itself 
given by the natural restriction maps
$$\iota_{st}:=e_{st}:E_t \to E_s,\;\; \textup{for}\;\;t>s $$
of $E$. Two local maps $u,v$ are called {\em inverse} to each other if
$$u \circ v=v\circ u=\iota $$
so that $\iota$ plays the role of the identity mapping.
%Note also that a Kolmogorov space morphism
%$$u:E \to F$$
%induces morphisms
%$$u_{st}:E_t \to F_s,\ x \mapsto f_{st} \circ u_t(x)$$
%and that this morphism is $0$-local if the norm of $u_t$ admits a uniform bound% in $t$.

\subsubsection{Exponential mapping}
\label{SSS::exponential}
Given a power series
$$ f=\sum_{n \geq 0} a_n z^n$$
we use the notation
$$ |f|=\sum_{n \geq 0} |a_n| z^n.$$

The {\em Borel transform}  of a formal power series is defined by
$$B:\CM[[z]] \to \CM[[z]],\ \sum_{n \geq 0} a_n z^n \mapsto \sum_{n \geq 0} \frac{a_n}{n!} z^n. $$
Consider a K1-space $E$ and its associated K2-space $ \Lt^1(E,E)$ of 1-local
operators. We define the set
\[ \Xt(R):=\{(t,s,u) \in \Lt^1(E,E)\;|\;\; ||u|| < R (t-s) \}\] 

\begin{theorem}\label{T::Borel} Let
$f=\sum_{n \geq 0} a_n z^n \in \CM\{z\}$ be a convergent power series with $R$ as radius 
of convergence. Then there is a well-defined map of spaces over $\D \subset \RM^2$, that 
we call the {\em Borel map}:
$$\Bt f:\Xt(R) \to \mathcal{H}{om}(E,E),\ (t,s,v) \mapsto (t,s,\sum_{n=0} \frac{a_n}{n!} v^n) $$
and one has the estimate
$$\|\Bt f(u)\| \leq |f|\left(\frac{||u||}{t-s}\right) .$$
\end{theorem}
\begin{proof}
For $u \in \Lt^1(E,E)$, we have $u^n \in \Lt^n(E,E)$ and 
\[ \|u^n\| \leq \|u\|^n \]
There is an inclusion 
$$\p_n: \Lt^n(E,E) \to \mathcal{H}{om}(E,E) $$
of Banach-spaces over $\Delta$, for which we have the bound
$$\frac{n^n}{e^n(t-s)^n}$$
on the norm. As $n^n \le e^n n!$, the choice of the constant $e$ leads to
a simple estimate for this norm:
$$\| \p_n\|=\frac{n^n}{e^n(t-s)^n} \leq \frac{n!}{(t-s)^n} .$$
We obtain:
\[  \|\sum_{n=0}^{\infty} \frac{a_n}{n!} \p_n(u^n)\| \leq \ \sum_{n=0}^{\infty} |a_n| \left(
\frac{||u||}{t-s}\right)^n= |f|\left(\frac{||u||}{t-s}\right) .\]
This proves the theorem.
\end{proof} 

\begin{corollary} 
\label{C::product}
Let $E$  be an $S$-Kolmogorov space and $(t:=t_0,t_1,t_2,\ldots) $ a 
decreasing sequence converging to $s>0$. For any sequence
$(t_n,u_n) \in L^1(E,E)$  such that
\begin{enumerate}[{\rm i)}] 
\item $\| u_n \| \leq  t_n-t_{n+1}$
\item $\s:=\sum_{n \geq 0} || u_n ||/(t_n-t_{n+1}) < +\infty$ 
\end{enumerate}
the sequence 
\[g_0=e^{u_0},\;\;g_1=e^{u_1}e^{u_0},\]
$$g_n:= e^{u_n}e^{u_{n-1}}\cdots e^{u_0}  $$ 
converges to an element $ g \in \Hom(E_t,E_s)$. 
Furthermore, we have the estimate:
$$ | g-\iota_{st} | < \frac{\s}{1-\s},  $$
where $\iota$ is the restriction map.
 \end{corollary}
\begin{proof}
The maps $e^{u_i}$ define elements 
of the Banach space $Hom(E_{t_i},E_{t_{i+1}})$ as long as $\|u_i\| \le t_i-t_{i+1}$, 
which holds by the first assumption. As a consequence
the compositions
\[e^{u_0},\;\;e^{u_1}e^{u_0},\;\;\ldots,e^{u_n},\ldots ,e^{u_1} e^{u_0},\ldots\]
are well defined. Furthermore, the Borel estimate gives
$$|e^{u_i}| \leq \frac{1}{1-\nu_i},\ \nu_i:=\| u_i\|/(t_i-t_{i+1}) .$$
As 
$$\frac{1}{1-x} \times \frac{1}{1-y} < \frac{1}{1-(x+y)} $$
for $x,y \in ]0,1[$, we get for the composition $e^{u_{i+1}} e^{u_i}$
$$|e^{u_{i+1}} e^{u_i}| \leq \frac{1}{1-(\nu_i+\nu_{i+1})} .$$
By a straighforward induction (and the fact that restrictions have norm
$\le 1$), we obtain the estimate
$$| g_n  | \leq \frac{1}{1-(\sum_{i=0}^n\nu_i)} .$$
Therefore
$$ | g_{n+1}-g_n  | \leq  \frac{| e^{u_{n+1}}-\iota |}{1-(\sum_{i=0}^n\nu_i)}  $$
Using again the Borel estimate 
$$ |e^{u_{n+1}}-1| \le \frac{\nu_{n+1}}{1-\nu_{n+1}} $$ 
we get 
$$| g_{n+1}-g_n  | \leq 
\frac{ \nu_{n+1}}{1-(\sum_{i=0}^{n+1}\nu_i)}  $$
From this it follows that the sequence $g_n$ converges in the Banach
space $Hom(E_t,E_s)$ with operator norm and that
$$|g-\iota|  \leq  \frac{\sum_{i \geq 0}\nu_i)}{1-(\sum_{i \geq 0}\nu_i)}= \frac{\s}{1-\s}  $$
\end{proof}
%%%%%%%%%%%%%%
\subsection{Arnold spaces}
In applications to iterations one often encounters Banach-spaces $E_{n,t}$ 
indexed by a discrete iteration variable $n \in \Nb:=\NM \cup \{\infty\}$ 
and continuous variables $t$ controlling the size of some neighborhood.
There are restriction mappings $E_{n,t} \to E_{m,t}$ and $E_{n,t} \to E_{n,s}$
for $n \le m$ and $t \ge s$, making commutative diagrams
\[
\begin{array}{ccc}
E_{n,t}&\longrightarrow&E_{n,s}\\
\downarrow&&\downarrow\\
E_{m,t}&\longrightarrow&E_{m,s} .\\
\end{array}
\]
Such a structure can be seen as a Kolmogorov space over a base of the form $\Nb \times B$, 
where of course we have to use the opposite ordering on the first variable. We call this structure
an {\em Arnold space}. It can also be seen as a (compatible) sequence $E_n$, indexed by $ n \in \Nb$ 
of ordinary Kolmogorov spaces, or better, as a functor $F:\Nb \to {\bf Kol}$, 
where ${\bf Kol}$ denotes the category of Kolmogorov spaces.
In particular there is fibre-wise notion of locality: a local map $u \in L^k(E,F)$ is a family of 
local maps
$$u_n \in L^k(E_n,F_n)$$ and thus defines a {\em norm sequence} $(|u_n|)$ which in applications needs to be controlled.
A descreasing sequence $n \mapsto s_n$ defines a map 
$$\s: \Nb \to \Nb \times B,\;\;\;n \mapsto (n,s_n),$$ 
that can be used to 'pull-back' an Arnold space $E$ and form a Kolmogorov space $E':=\s^*E$ over
$\Nb $ with fibre $E_{n,s_n}$. 
The consideration of Arnold-spaces is useful in situations where one wants to postpone 
the choice of an appropriate sequence $(s_n)$ as long as possible.\\

In this paper we consider Banach spaces of holomorphic functions $\Ot^k$ on sets $W_{n,s}$.
These sets combine into a relative open set $ W \to \Nb \times B$ and the Banach spaces
combine into an Arnold space $\Ot^k(W)$.

\begin{proposition}Let $u \in L^m(E,F)$, $v \in L^n(F,G)$ be local maps with norm subquadratic norm 
sequences then the norm sequence $u \circ v \in  L^{m+n}(E,G)$ is subquadratic with order at most
$$\ord (|u \circ v |) \leq \max(\ord(|u_n|),\ord(|v_n|)). $$
\end{proposition}
The proof is obvious.

%%%%%%%%%%%%%%%%%%%%%%%%%%%%%%%%%%%%%%%%%%%%%%%%%%%%%%%%%%%%%%%%%%%%%%%%%%%
\section{Kolmogorov spaces in the analytic context}
\subsection{Spaces of holomorphic functions}
\subsubsection{}
Let $U$ be a relative compact open subset of $\CM^n$, and $\Ot(U)$ the ring
of holomorphic functions on $U$. It is naturally a Fr\'echet space.
One can also attach to $U$ various Banach spaces of holomorphic funcions.
The space space of {\em bounded holomorphic functions} is denoted by
\[\Ot^b(U):=\{f \in \Ot(U)\;|\;\;f\;\textup{is bounded} \}\]
is a Banach space with 
\[ |f|=\sup_{z \in U}|f(z)|\]
as norm. The space of {\em square integrable holomorphic functions}, 
denoted by
\[\Ot^h(U):=\{f \in \Ot(U) \;|\;\; \int_{\overline{U}} |f|^2 dV \}, \]
is a Hilbert space with the $L^2$-norm as norm.
We denote by 
\[\Ot^{k}(U)\] 
the Banach space of complex valued (Whitney) {\em $C^k$-functions on the closure of
$\overline{U}$}, which are holomorphic on the interior of $U$ with 
\[|f|=\max_{|I| \leq k}\sup_{z \in U}|\d^I f(z)|\]
as norm.\footnote{Note that $\Ot^0(U)=:\Ot^c(U)$ is the same as the space of holomorphic 
functions that extend continuously to the boundary.}
As any $C^k$-function is bounded and any bounded function on a relative compact set is square-integrable, there are natural inclusions
\[ \Ot^k(U) \subset \Ot^b(U) \subset \Ot^h(U).\]

If we generalise this to the context of relative open sets over some base, 
we obtain the most important class of Kolmogorov spaces.

\subsection{}
Let $\St$ denote the set of subsets of $\CM^n$, partially 
ordered by inclusion.

\begin{definition} Let $(B,\ge )$ be a partially ordered set.
A {\em set over $B$} in $\CM^n$ is an order reversing map 
$$B \to \St$$
\end{definition}

So a set over $B$ consists of sets $U_t $ for $t \in B$ such that for
$t \ge s$ we have an inclusion $U_s \hookrightarrow U_t$.
(It can be seen as a contravariant functor if we consider $B$ and $\St$ 
as categories.) If all sets $U_t$ are open/closed, we call it an {\em open/closed 
set over $B$}.  We will write such an open set over $B$ as $U \to B$, with 
fibres $U_t$ over $t \in B$.

As a simple example, the {\em relative unit polydisc} $D \to \RM_{>0}$ with fibres
$$D_t=\{ z \in \CM^d: |z_1| \leq t, \dots,|z_{d}| \leq t\}. $$

Given sets $U \to B, \ U' \to B$ we may perform many of the
usual operations fibre-wise. For example we may form their fibred product
$$U \times_B U' \to B , $$
with fibres the cartesian product of the fibres of $U$ and $U'$, etc.

%%%%%%%%%%%%%%%%%%%%%%%%%%%%%%%%%%%%%%%%%%%%%%%%%%%%%%%%%%%%%%%%%%%%%%%%%%%
\subsubsection{}
 Using such relative open sets we can create a plethora of Kolmogorov spaces.
For an open set $U \to B$ over $B$ we may for each $t \in B$ form the 
Banach space $ \Ot^b(U_t)$ of {\em bounded holomorphic functions} on $U$, 
with the sup-norm as norm. There are for $s \le t$ obvious restriction maps 
$\Ot^b(U_t) \to \Ot^b(U_s)$ of norm $<1$, 
and hence we obtain a Kolmogorov-space
$$\Ot^b(U) \to B $$
with $\Ot^b(U)_t:=\Ot(U_t)$. Similarly, the spaces $\Ot^h(U_t)$ of square 
integrable and $\Ot^k(U_t)$ of $C^k$-functions for Kolmogorov spaces $\Ot^h(U)$ 
and $\Ot^k(U)$. There are natural Kolmogorov space morphism:
$$ \Ot^k(U) \to \Ot^b(U), \;\;\;\Ot^b(U) \to \Ot^h(U) .$$

\subsubsection{Cauchy-Nagumo estimate.}
Let $D_t \subset \CM$ be the disc of radius $t$. For a holomorphic function 
$f \in \Ot^c(U_t)$ one has the following elementary estimate
\[ |f^{(m)}|_s \le \frac{m!}{(t-s)^k}|f|_t\]
for $s < t$, which is a straightforward consequence of the Cauchy 
integral formula and differentiation under the integral sign. 

This simple idea can be extended to general partial differential operators
on appropriate relative open sets.

\begin{definition} We say that an open and relatively compact 
set $U$ of $\CM^n$ over $B=]0,S]$ is a {\em Huygens set}, if 
the following condition holds
$$\forall x \in U_s,\ x+D_{(t-s)} \subset U_t,\ a \in \RM $$
for any $s<t \le S$.\\
\end{definition}

The proof of the 1-variable case has an immediate generalisation to

\begin{proposition} ({Cauchy-Nagumo})
Let $B=]0,S]$ and $U$ a Huygens set over $B$, then any partial differential 
operator 
$$P=\sum_{|I| \le m} a_I(z)\d^I, a_I \in \Ot^k(U_S)$$
of order $m$ defines an $m$-local operator of the Kolmogorov space $\Ot^k(U)$
over $B$:
\[ P \in L^m(\Ot^k(U),\Ot^k(U))\]
\end{proposition}

In applications one encounters often slightly more general situations, like
the following.

\begin{definition} Let $(a_n) \in \RM_{>0}$ be a falling positive sequence. We say that an open and relatively 
compact set $U$ of $\CM^n$ over $B=]0,S] \times \Nb $ is an {\em $a$-Huygens set}, if 
the following condition holds
$$\forall x \in U_{n,s},\ x+D_{a_n(t-s)} \subset U_{n,t}$$
for any $s<t \le S$ and $n \in \NM$.
\end{definition}
The sequence may very well
go to $0$, making an uniform choice for $a$ impossible.

The above proposition admits a straightforward variant:
\begin{proposition} ({Cauchy-Nagumo II})
\label{P::Cauchy_Nagumo}
Let $a=(a_n)$ be a falling subquadratic sequence.
Let $B=]0,S] \times \Nb$ and $U$ an $(a_n)$-Huygens set over $B$, then any partial differential 
operator 
$$P=\sum_{|I| \le m} a_I(z)\d^I, a_I \in \Ot^k(U_S)$$
of order $m$ defines an $m$-local operator of the Kolmogorov space whose norm sequence $|P|_n$ is subquadratic with index bounded by that of $(a_n)$:
$$\ord(|P|_n) \leq \ord(a) .$$
\end{proposition}

\subsubsection{}
\label{SSS::local_equivalence}
A function $f \in \Ot^h(U)_t$ is by definition holomorphic in $U_t$. If $U$ is
Huygens, then given any $s<t$, its restriction to $\overline{U}_s$ is a 
$C^k$-function for any $k$. So we have 'restriction mappings' 
\[J_{st}:\Ot^h(U_t) \to \Ot^k(U_s)\]
These maps are compatible with the restrictions on $\Ot^h(U)$ and $\Ot^k(U)$,
so combine into an element
\[ J \in Hom_\D(\Ot^h(U),\Ot^k(U))\]
 
\begin{proposition}
\label{P::local_equivalence}
If $U$ is an $a=(a_n)$-Huygens set over $B$, then 
$$J: \Ot^h(U) \to \Ot^k(U) $$
is a local map. Moreover, if $a$ is a falling subquadratic sequence, the norm sequence $|J_n|$ is bounded by an increasing
subquadratic sequence with the same exponent.
\end{proposition}
\begin{proof}
Consider a function $f \in \Ot^h(U)_t$ and let $s <t$. The Taylor expansion 
of $f$ at a point $w \in U_{n,s}$ reads:
$$f(z)=\sum_{J \in \NM^d} a_J (z-w)^J,\ a_J \in \CM, $$
by assumption,  $U$ is an $(a_n)$-Huygens set, so the polydisc $D_w$ centred at $w$ 
with radius $\s= a_n(t-s)$ is contained in $U_t$. We then have
$$\int_{D_w} | f(z)|^2 dV=\sum_{J \in \NM^d} C(J) |a_J|^2 \s^{2|J|+2d},\;\;C(J)=\prod_{k=1}^d \frac{\pi}{j_k+1}, $$
where $dV$ is the Lebesgue measure.\\

So we obtain
\[C(0) |a_0|^2 \s^{2d} \le \int_{D_w} | f(z)|^2 dV \leq  \int_{U_t} | f(z)|^2 dV=| f |_t^2 . \]
This shows that
$$ |f(w) |=|a_0| \leq  \frac{c}{\s^d}\left( \int_{\D_w} | f(z)|^2 dV \right)^{1/2} \leq   \frac{c}{a_n^d(t-s)^d} | f |_t$$
for any $w \in U_{n,s}$ and $c:=\sqrt{\frac{1}{C(0)}}$.
When we apply the same argument to the derivatives of $f$ and combine it
with the Cauchy-Nagumo estimate, we find an estimate of the form
\[ | Jf |_s =\max_{|I| \le k }\sup_{w \in U_s}|\partial^If(w)| \leq \frac{c'}{a_n^{d+k}(t-s)^{d+k}} | f |_t.\]
 The proposition follows.
\end{proof}
%%%%%%%%%%%%%%%%%%%%%%%%%%%%%%%%%%%%%%%%%%%%%%%%%%%%%%%%%%%%%%%%%% 
\subsection{Kolmogorov spaces attached to the torus $\TM$}

\subsubsection{}
We denote by
\[T_t:=\{ z \in (\CM^*)^n:\ e^{-t}<| z_i| \leq e^t\} \]
the neighbourhood of $\TM \subset \CM^{n}$. 
The space of holomorphic functions $\Ot(T_t)$ can be identified with the
analytic Fourier series
\[ f=\sum a_I z^I \in \CM\{z,z^{-1}\}\]
for which
\[ |a_I| = O(e^{-|I| t}).\]
The sets $T_t$ can be seen as fibres of an open set $T$ over $\RM_{>0}$. When we
restrict it to $B=]0,s_0]$, it is an $a$-Huygens set, for an appropriate $a$
(which goes to $0$ if $s_0 \to \infty$). We will consider the corresponding 
Kolmogorov spaces
\[\Ot^k(T),\;\;\;\Ot^b(T),\;\;\;\Ot^h(T)\]
over $]0,s_0]$. Clearly, the elements of each of the underlying Banach spaces
$\Ot^h(T_t)$ can be seen a special elements of $\CM\{z,z^{-1}\}$.

\subsubsection{The Arnold-Moser lemma}
\label{SSS::Arnold_Moser}
\begin{lemma} Assume that a function $f \in \Ot^h(T)_s$ depends only on 
harmonics of degree $ \ge m$, then for $s \le t$ we have the estimate
$$| f|_s \leq \left(\frac{e^s}{e^{t}} \right)^{m/2} | f |_t $$
\end{lemma}
\begin{proof}
The cartesian product
\[S_t:=S_t(1) \times S_t(2) \times S_t(2d)\]
of coordinate strips
\[ S_t(j):=\{ x_j=\xi+i\eta_j\;|\;\;0<\xi_j\le 2\pi,\;\;\;-t <\eta_j <t\}\]
parametrises the torus neighborhood $T_t$ via the maps 
\[x_j \mapsto z_j=e^{i x_j} .\]
We use the $L^2$-norm on $\Ot^h(T_t)$, obtained by integration of the pull-back of $f(z)\overline{f(z)}$ over the strip $S_t$. The monomials $z^I$ then 
form an orthogonal basis. As in one variable we have 
\[\int_{S_r} e^{inx} \overline{e^{inx}}d\xi d\eta =\int_{S_r} e^{-2n\eta} d\xi d\eta = \left\{ \begin{array}{l} 2\pi \cdot 2r\;\;\text{if}\;\;n=0\\2\pi \frac{\sinh(2nr)}{n}\;\;\text{if}\;\;n \neq 0 \end{array} \right.\]

we find that for $I=(i_1,i_2,\ldots,i_{2d})$
$$| z^I |_t^2=(2\pi)^{2d}\prod_{k=1}^{2d} \frac{\sinh(2i_k t)}{i_k}$$
By the Pythagorean theorem, for $f \in \Ot^h(T_t)$, we have:
\begin{align*}
| f |_s^2&= \sum_{|I| \geq m} |a_I|^2|z^I|_s^2  \\
             &= \sum_{|I| \geq m} |a_I| \frac{|z^I|_s}{|z^I|_t} |z^I|_t\\
             & \leq  \frac{\sinh(2m s)}{\sinh(2m t)}|f|_t^2 .
\end{align*}
Here we used the two inequalities
\[  \frac{|z^I|_s}{|z^I|_t} \le \frac{\sinh(2|I|s)}{\sinh(2|I|t)} \le \frac{\sinh(2 m s)}{\sinh(2 m t)}.\]
The first on is implied by the fact that for fixed positive numbers 
$a,b,\ldots,z>0$ the function
\[ x \mapsto \frac{\sinh(ax)\sinh(bx)\ldots\sinh(zx)}{\sinh((a+b+\ldots+z)x)}\]
is monotonous increasing in $x$. The second inequality follows because all
$I$ appearing in the sum have $|I|\ge m$ and the function
\[ x \mapsto \frac{\sinh(ax)}{\sinh(bx)}\]
is monotonous increasing in $x$ for $a \ge b >0$.
Finally, as $t>s$, one has also: 
$$  \frac{\sinh 2 m s}{\sinh 2 m t} =\frac{e^{ms}(1-e^{-2ms})}{e^{mt}(1-e^{-2mt})} \leq \frac{e^{ms}}{e^{mt}}$$
\end{proof}

\subsubsection{Differential forms }
We can define similarly Kolmogorov spaces of relative one-forms
by putting 
\[\Omega^{k,1}(T):=\Ot^k(T) \frac{dz_1}{z_1}\oplus \Ot^k(T) \frac{dz_2}{z_2}\oplus\ldots \oplus \Ot^k(T) \frac{dz_n}{z_n}\]
where we define the norm of a form $\a=\sum a_i \frac{dz_i}{z_i}$ to be
\[ |\a |_t:= \sup_{1 \le i \le n} \{|a_i|_t \} .\]
By setting
$$\Omega^{k,l}(T):=\wedge^l \Omega^{k,1}(T)$$
it can be extended to higher values of $l$, and there are similarly versions
for $b, h$ instead of $k$.

\subsubsection{De Rham Complex}
By \ref{P::Cauchy_Nagumo}, the exterior derivative 
$$d:\Omega^{k,l}(T) \to \Omega^{k,l+1}(T) $$
is a 1-local morphism. From this it follows that the space of
closed forms 
$$Z^{k,l}(T) \subset \Omega^{k,l}(T) $$ 
form a Kolmogorov subspace. 

The $1$-forms $\a_k=\iota_{\theta_k}\omega$ are de Rham dual to 1-cycles  
$$\g_1,\dots,\g_{2d} \in H_1(T,\CM)=H_1(\TM,\CM). $$
On $Z^{k,1}$ we define linear forms
\[ c_k: Z^{k,1}(T) \to \CM,\ \a \mapsto \int_{\g_k} \a .\]
The subspace $B^{k,1}(T)\subset Z^{k,1}(T)$ of exact 1-forms
coincides with the forms with vanishing period integrals,
so  the $1$-form
$$\b=\a-\sum_{k=1}^{2d} c_k(\a)\a_k$$
belongs to the space $B^{k,1}$.
Consider the path 
$\g_z$ connecting $|z|:=(|z_1|,|z_1|,\ldots,|z_n|)$ to $z \in T_s$ 
by changing only the arguments: 
\[ \g_z:[0,1] \to T_s,\;\;t \mapsto (|z_1|e^{ i t \theta_1},\ldots,|z_n|e^{i t \theta_n} )\]
where $\theta_i=arg(z_i) \in ]0,2\pi]$.
Integration over $\gamma_z$ defines a map of Kolmogorov spaces
$$\int: B^{k,1}(T) \to \Ot^k(T),\ \b \mapsto [z \mapsto \int_{\g_z} \b],$$
as the exactness of the form $\b$ guarantees that the function $\int \b$ 
is continuous. The obvious estimate
$$|\int_{\g_z} \b| \leq \left(2\pi e^s\right)^{2d}|\b| \leq \left(2\pi e^{s_0}\right)^{2d}|\b|$$
guarantees boundedess of the map.\\

We have shown the
\begin{proposition}
The maps $c_k$, $\int$ define a morphism of Kolmogorov spaces:
\begin{align*}
  Z^{k,1}(T) & \to  \CM^{2d} \oplus \Ot^k(T) ,\\
\a &\mapsto (c_1(\a),\dots,c_{2d}(\a),\int \b)
 \end{align*}
where 
\[\beta:=\a-\sum_{i=1}^{2d}c_i(\a)\a_i\] 
\end{proposition}

\subsubsection{Symplectic vector fields}
\label{SSS::symplectic}
 The Kolmogorov space of relative vector fields is defined to be dual 
to $\Omega^{k,1}$:
\[\Theta^k(T):=\Ot^k(T)\theta_1 \oplus \Ot^k(T)\theta_2\oplus \ldots \oplus \Ot^k(T) \theta_n,\]
with $\theta_j:=z_j\d_{z_j}$
By \ref{P::Cauchy_Nagumo}, there is an embedding of Kolmogorov spaces
$$\xymatrix{\Theta^k(T)  \ar@{^{(}->}[r] & L^1(\Ot^k(T),\Ot^k(T))}  .$$

As the Lie bracket is a first order differential operator in the coefficients, 
\ref{P::Cauchy_Nagumo} also implies a functorial analytic version of the adjoint representation:
\begin{proposition} The Lie bracket defines a 1-local map
$$ad:\Theta^k(T) \to L^1(\Theta^k(T), \Theta^k(T)),\ X \mapsto [X,-] .$$
\end{proposition} 
% 
%  \begin{proposition}
%  The map 
% $$ L^1(\Ot^b(T),\Ot^b(U)) \to  \Theta^b(U),\ v \mapsto (z_1^{-1}v(z_1),\dots, z_{2d}^{-1}v(z_{2d}))$$
% is a bounded  right inverse of the above one.
% \end{proposition}
% \begin{proof}
% Put
% $$v(z)=\sum_{k \geq 0} v_k(z)\d_{z_k}$$
% Our task is to show that each component $v_k$ is bounded in $U_t$. Choose $r<s<t$ so that
% $s$ is the midpoint of $[r,t]$. The maximum of $|v_k|$ in $U_r$ is reached at a point 
% $$m=(re^{i\theta_1},\dots,re^{i\theta_{2d}})$$
% We consider the function
% $$f_k(z):=\frac{1}{e^{i\theta_k}t-z_k} \in \Ot^b(U)_s. $$
% As
% $$\d_{z_k}f_k(z)=\frac{1}{(e^{i\theta_k}t-z_k)^2}  $$
% We have:
% \begin{align*}
% v(\frac{1}{e^{i\theta_k}t-z_k})&=\sum_{k \geq 0} v_k(z)\d_{z_k}\frac{1}{e^{i\theta_k}t-z_k}\\
% &=z_k v_k(z)\frac{1}{(e^{i\theta_k}t-z_k)^2} 
% \end{align*} 
% and estimating the right-hand side at $z=m$ we get that
% \begin{align*}
% \left|v(\frac{1}{e^{i\theta_k}t-z_k})\right|_r& \geq |m_k v_k(m)\frac{1}{(e^{i\theta_k}t-m_k)^2}|_r\\ 
%   &= \frac{ (1-r)| v_k(m)|}{(t-r)^2}
% \end{align*}
% As $v$ is $1$-local we also have
% \begin{align*}
% |v(\frac{1}{e^{i\theta_k}t-z_k})|_r & \leq \frac{e|v|_t}{(s-r)|e^{i\theta_k}t-z_k|_s}=\frac{4e|v|_t}{(t-r)^2}
% \end{align*}
% Comparing these two estimates we get that
% $$|v_k|_r=|v_k(m)| \leq \frac{4e|v|_t}{1-r}\leq \frac{4e|v|_t}{1-b}.  $$
% This proves the proposition.
% \end{proof}

By the standard symplectic duality isomorphism,
\[ \Theta \to \Omega^1,\;\;\; X \mapsto \iota_X(\omega) \]
the closed forms correspond to symplectic vector fields
and the exact forms to Hamiltonian fields.

The decomposition into exact and non-exact part of a closed $1$-form 
translates into the following statement

\begin{proposition}
The decomposition of symplectic derivations into exact and non-exact part:
\[\St^k(T) \to \Ot^k(T) \oplus \bigoplus_{j=1}^{2d} \CM \theta_j ,\]
$$X \mapsto \int \iota_X\omega + \sum_{j=1}^{2d}c_j(\iota_X\omega) \theta_j$$
is a locally bounded morphism of Kolmogorov spaces.
\end{proposition}
%%%%%%%%%%%%%%%%%%%%%%%%%%%%%%%%%%%%%%%%%%%%%%%%%%%%%%%%%%%%%%%%%%%%%%%%
%%%%%%%%%%%%%%%%%%%%%%%%%%%%%%%%%%%%%%%%%%%%%%%%%%%%%%%%%%%%%%%%
\subsection{Functional spaces involved in the iteration}
\subsubsection{Frequencies}
 We use
the following local variant of arithmetic classes: 

\begin{definition} For a fixed $\nu$ and falling sequence $a=(a_n)$ and
$s_0 \in \RM_{>0}$ we define a sets
\[Z_{n,s}:=Z_{n,s}(\nu,a,s_0):=\{\phi \in D_s^{2d}: \forall k \leq n,\s(\nu+\phi)_k  \geq  a_k(s_0-s) \} \]
\end{definition}
It is readily checked that for $n \le m$ and $s \le t \le s_0$ one has:
\[ Z_{m,s} \subset Z_{n,s},\;\;\; Z_{n,s} \subset Z_{n,t},\]
so that the $Z_{n,s}$ are fibres over a relative set
$$Z(a) \to \Nb \times ]0,s_0] .$$

\begin{lemma} \label{L::calibrated} The set $Z(a)$ is an $a^*$-Huygens set:
$$Z_{n,s}+ D_{a^*_n(t-s)}  \subset Z_{n,t} ,$$
where
\[ a^*_n:=\frac{a_n}{2^n}\]
\end{lemma}
\begin{proof}
Let  $\phi \in Z_{n,s}$ and take $x \in \CM^{2d}$  satisfying 
$$\| x\|  \leq \frac{a_n}{2^n} (t-s).$$ 
For any $k \le n$ and $\| J \| \le 2^k$ we then have:
\[|(x,J)| \le \|x\| \|J\| \le \frac{a_n}{2^n} (t-s)\cdot 2^n= a_n(t-s)\le a_k(t-s).\]
So we obtain
\begin{align*}
|(\nu+\phi+x,J)|     & \geq |(\nu+\phi,J)|-|(x,J)| \\
                     & \geq a_k(s_0-s)- a_k(t-s)=a_k(s_0-t) .\\                          
\end{align*}
This shows that $\phi+x \in Z_{n,t}$ and thus proves the lemma.
\end{proof}
\begin{corollary} Assume $a=(a_n)$ is a falling subquadratic sequence then any partial differential operator with coefficients in
$\Ot^c(Z)$ defines a local map whose norm sequence is bounded by a rising subquadratic sequence with the same index.
\end{corollary}

%%%%%%%%%%%%%%%%%%%%%%%%%%%%%%%%%%%%%%%%%%%%%%%%%%%%%%%%%%%%%%%%
\subsubsection{Functional spaces with parameters}
We have considered the relative neighbourhood of $\TM$
$$T \to \RM_{>0}$$
with fibre
$$T_s=\{ z \in (\CM^*)^{2d}:\forall i,\ e^{-s}< |z_i|< e^s\} $$
For the formulation of the iteration we introduced detuning variables
\[\phi_1,\phi_2,\ldots,\phi_{2d}\]
and variables 
\[\dt_1,\dt_2,\ldots,\dt_{l}\]
describing the perturbation of the symplectic form. We have to introduce
appropriate  neighbourhoods in the space with variables $\phi,\dt, z$.

\begin{definition}  For  fixed decreasing sequence $a$ and $s_0$ we set
$$W(a) \to \Nb \times \RM_{>0},\  V \to \Nb \times \RM_{>0}$$
by putting \[W(a):= Z(a) \times_{\RM_{>0}} D^l \times_{\RM_{>0}} T,\;\;\; 
V(a):=Z(a) \times_{\RM_{>0}} D^{l} \]
where $D \to \RM_{>0}$ is the relative unit polydisc. 
\end{definition}

The coordinates on $W_{n,t}$ are 
\[\phi_1,\ldots,\phi_{2d},\;\dt_1,\ldots,\dt_{l},\;z_1,\ldots,z_{2d} ,\]
and there are projection maps 
\[ W(a) \to V,\;\;\;(\phi,\dt,z) \mapsto (\phi,\dt).\]
The functional spaces we consider are the Arnold spaces 
$$\Ot^k(W(a)),\Ot^h(W(a)),\;\; \Ot^k(V),\ \Ot^h(V),\;\; \Theta^k(W(a)),\ \St^k(W(a)),\ etc.$$
over $\Nb \times \RM_{>0}$.

%%%%%%%%%%%%%%%%%%%%%%%%%%%%%%%%%%%%%%%%%%%%%%%%%%%%%%%%%%%%%%%%%%%%%%%%%%%%%%%%%%%
\subsubsection{Model iteration}\label{SSS::model}
Bruno sequences appear naturally in connection with quadratic iterations
of the type
\[ x_{n+1} =a_n x_n^2 .\]
As solution one has
\[x_1=a_0 x_0,\;\;x_2=a_1x_1^2=a_1a_0^2x_0^2=(a_0a_1^{1/2} x_0)^2\]
and is solved in general by
\[ x_n=(a_0a_1^{1/2}\ldots a_n^{1/2^n} x_0)^{2^n},\]
so that the sequence $(x_n)$ converges quadratically to $0$ if $x_0 < 1/a_{\Pi}$.

Our aim is to investigate a slightly more general iteration of the form
$$x_{n+1}=\frac{1}{2}\left( a_n x_n^2+b_n x_n\right) $$
\begin{proposition}
\label{P::model}
Let $a=(a_n) , a_n \ge 1$ be an increasing and $b =(b_n)$ a falling sequence
of poritive numbers. Assume that for some $N \in \NM$ one has
$$(\star)_N\qquad n \geq N \implies b_n^2a_n \leq b_{n+1}. $$
Then for any $0 <\e$ there exists $\dt$ such that $x_0 \leq \dt$ the real sequence
$$x_{n+1}=\frac{1}{2} \left(a_nx_n^2+ b_n x_n \right)$$
converges to zero and one has the estimate
$$x_{n} \leq \e\, b_n .$$
\end{proposition}
\begin{proof}
Without loss of generality we may assume that $\e \leq 1$. We will give
an explicit $\dt$ that does the job. Let us start by remarking if we
have $x_n \le \e b_n$ for some $n \ge N$,  
then we find, using $(\star)_N$, $a_n \ge 1$ and $\e^2 \le \e$ 
\begin{align*}
x_{n+1}= \frac{1}{2} \left(a_nx_n^2+b_n x_n \right)& \leq   \frac{1}{2} \left(a_n\e^2\,b_n^2 +\e\,b_n^2  \right) \leq \e\, a_nb_n^2 \leq   \e\,b_{n+1}
\end{align*}
so that if $x_N \le \e b_N$ holds, then it holds for all $n \ge N$.
We will define constants $c_0,c_2,\ldots, c_N  \le 1$
such that so that if 
\[ x_0 \le c_0 \cdot c_1 \cdots c_2 \cdots c_N\]
holds, then for $n=1,2,\ldots,N$ we have 
\[ x_n \le c_n \cdot c_{n+1} \cdots c_{n+2} \cdots c_N .\]
Assuming this holds true at rank $n<N$, then
\begin{align*}
 a_nx_n^2 & \le a_n c_n^2c_{n+1}^2\ldots c_N^2 \le (a_n c_n) c_{n+1}\ldots c_N\\
 b_n x_n & \le (b_n c_n) c_{n+1}\ldots c_N \\        
x_{n+1}&\le (\frac{a_n+b_n}{2}c_n) c_{n+1} \ldots c_N 
\end{align*}
So if we take
\[c_n \le \min(2/(a_n+b_n),1)\]
we find $x_{n+1} \le c_{n+1}\ldots c_N$. So if we put finally 
\[c_N \le \min(\e b_N, 1)\]
then the number 
\[ \dt:=c_0\cdot c_1 \cdot \ldots c_{N}\]
does the job. This concludes the proof of the proposition.
\end{proof}

We will use the above Proposition for subquadratic sequences $a=(a_n)$ and $b=(b_n)$ of the form
$$a_n=Ae^{B\a^n},\;\;\; b_n=Ce^{-D\b^n},\ \b>\a \text{ and } \b \in ]1,2[ $$ 
and where $A,B,C,D \in \RM_{>0}$. Condition $(\star)_N$ is obviously fullfilled.
 
\section{Proof of the Stability Theorem}
\subsection{Final preparations}
\subsubsection{Definition of $j$ and $\sigma$}
The iteration scheme in the ring $R$ for the normal form of $\S 3$ 
was based on maps
\[j_n: \St(R) \to \Qt(R),\]
which were defined in terms of decomposition and truncation. It is not
difficult to lift these maps to the functional analytic level.

We fix a sequence $a$ and consider the Arnold-space $\St^k(W(a))$. Its components $\St^k(W(a))_n$ 
are K1-spaces of symplectic vector fields $X$ on sets $W_{n,t}$.
For these we first  make a {\em decomposition} as in \ref{SSS::symplectic}:
$$ X =\{-,f\} +\sum_{i=1}c_i \theta_i, $$
with 
\[ f \in \Ot^k(W_{n,t}),\;\;\;c_i \in \Ot^k(V_{n,t})\]
and define the $A$-trunction as before:
\[ [X]_A :=\{-,[f]_A\}+\sum_{i=1}^{2d} [c_i]_A \theta_i \in \St^k(W_{n,s}).\]
where $s <t$. In the iteration we will use maps
\[\s(n,t,s):\St^k(W_{n,t}) \to \St^k(W_{n,s})\]
\[j(n,t,s): \St^k(W_{n,t}) \to ,L^1(\Theta(W),\Theta(W))_{n,s} \]
by setting	:
\[\s(n,t,s)(X):=[X]^{2^n}=\{-,[f]^{2^n}\}+\sum_{i=1}^{2d}[c_i]^{2^n}\theta_i\]
\[j(n,t,s)(X):=ad(\{-,[f]^{2^n} \star H\} +\sum_{i=1}[c_i]^{2^n}\d_{\phi_i}).\]
where $H$ is the resolvente and $ad$ is the adjoint action, which embeds the space $\Qt(W)$ inside $L^1(\Theta(W),\Theta(W))$ via a local map (see \ref{SSS::symplectic}).

It is easy to see that norm sequence of $|j_n|$ is subquadratic and satisfies
$$\ord(|j_n|) \leq \ord(a). $$
By the Cauchy-Nagumo lemma (see \ref{P::Cauchy_Nagumo}), the norm sequences $|\d_{\phi_i}|$ as well as the non-exact term have the same property:
 $$\ord(|\d_{\phi_i}|) \leq \ord(a). $$
Due to local equivalence taking the convolution with $H$ has again the same property.  The Poisson bracket being a biderivation, due to the Cauchy-Nagumo lemma, the exact term it is a local operator with bounded norm.

These maps form compatible systems and combine to form local maps (we take the same index for locality $p$ to simplify further estimates):
\vskip0.1cm
\begin{align*} \s & \in L^p(\St^k(W(a)),\St^k(W(a))) \\
 j& \in L^p(\St^k(W(a)),L^1(\Theta^k(W(a)),\Theta^k(W(a))))
 \end{align*}
If the sequence $a$ that goes in the definition of $W(a)$ is subquadratic  then $W(a)$ is an $a^*$-Huygens set,
the norm sequence of the map $j$ is subquadratic with the index smaller than that of $a$ , as it is the composition of such maps.  

%We have an estimate of the form:
%\[ | j(n,t,s)| \le  \frac{e^p\|j_n\| }{a_n^p (t-s)^p}=\frac{Ce^p }{a_n^{p+1} (t%-s)^p}\]
%for some constant $C$ (the $+1$ in the exponent of $a_n$ comes from the  divisi%on map).\\
%Due to the Arnold-Moser lemma \ref{SSS::Arnold_Moser}, the truncation map $\s$
%is not only tamed but is also 'small', as it admits an estimate of 
%the form
%$$| \s(n,t,s) | \leq \left(\frac{s}{t} \right)^{2^n} \frac{C}{a^p_n(t-s)^p} $$
%for $s< t \leq 1$.
%%%%%%%%%%%%%%%%%%%%%%%%%%%%%%%%%%%%%%%%%%%%%%%%%%%%%%%%%%%%%%%%%%%%%%%%
%%%%%%%%%%%%%%%%%%%%%%%%%%%%%%%%%%%%%%%%%%%%%%%%%%%%%%%%%%%%%%%%%%%%%%
\subsection{Pulling back}
We choose $\rho =(e^{-\b^n}) \in \SM^-$ with order  $\b \in ]1,2[$ higher than that of $(a_n)$ and therefore of  $|j_n|$. Then we define a sequence $(s_n)$ indexed by half-integer by putting
$$s_{n+1/2}:=\rho_n^{1/2^{n+1}}s_n, \;\;s_0=t $$
As $\rho$ is subquadratic, the sequence $(s_n)$ converges to a positive limit $s_\infty(\rho)$ and 
$$ s_n-s_{n+1/2}\sim \frac{1}{2} \left(\frac{\b}{2}\right)^n s_n \sim \frac{1}{2} \left(\frac{\b}{2}\right)^n s_\infty $$
 
 Pulling-back an Arnold space
$$E \to \Nb \times \RM_{>0}  $$
via the map
$$i:\frac{1}{2}\Nb \to \Nb \times \RM_{>0},\ n \mapsto (\lfloor n \rfloor,s_n)   $$
defines a Kolmogorov space that we denote by $\rho^*E$ over $\frac{1}{2}\Nb$.
Pulling back the maps $j$ and $\s$, we obtain maps whose norms sequence belong respectively to $\SM^+$ and $\SM^-$.
From the Arnold-Moser lemma (\ref{SSS::Arnold_Moser}), we deduce the estimate:
$$| \s(n,t,s) | \leq \left(\frac{e^{s}}{e^{t}}\right)^{2^{n-1}}\frac{C}{a^p_n(t-s)^p} $$
we deduce that:
$$| \s(n,s_{n},s_{n+1/2}) | \leq \left(\frac{e^{s_{n+1/2}}}{e^{s_n}}\right)^{2^{n-1}}\frac{C}{a^p_n(s_n-s_{n+1/2})^p} \sim \frac{C 2^{(n+1)p} e^{-\b^ns_\infty/4}}{a^p_ns_\infty^p\b_n^p} $$
Therefore the norm sequence $| \s(n,s_{n},s_{n+1/2}) |$ is a falling subquadratic sequence  of order~$\b$.
%%%%%%%%%%%%%%%%%%%%
\subsection{The Convergence Theorem}
\begin{theorem}
\label{T::normal_form}
Consider a symplectic vector field on the complex torus
\[ \Vt=\sum_{i=1}^n \nu_i z_i\d_{z_i}\]
and assume that $\s(\nu) \in \SM^-(\a), \a \in ]1,2[$. Consider a subquadratic sequences $a \in \SM^-(\a)$ such that
$a \leq \s(\nu)$ and define
$\rho=(e^{-\b^n})$ with $\b \in ]\a,2[$. Then for any $k,\e>0$, there exists $\dt>0$ such that for $ |S_0| \leq \dt$:
there is an iteration in $\rho^*\St^k(W(a))$ defined by
$$  S_{n+1}=f_*(j_n(S_n))\Vt+e^{- [j_n(S_n),-]}\s_n(S_n)$$
and it satisfies the estimate  $| S_n|<\e e^{-\b^n} $.
\end{theorem}
\begin{proof}
The analytic series
$$f(z)=e^{-z}(1+z)-1 \in \CM\{z\}$$
is the Borel transform of 
\[- \frac{z^2}{(1+z)^2} \in z^2 \CM\{z\} .\]
which has radius of convergence equal to $1$ and, choosing $r=1/2 <1$, we get that:
$$ \left|\frac{z^2}{(1+z)^2}\right| \leq \frac{|z|^2}{1-r^2}<2|z|^2$$
for $|z| \leq r$. Therefore assuming
$$ |S_n| \leq \frac{s_{n+1/2}-s_{n+1}}{2\|j_n\| } $$
 by \ref{SSS::exponential}, we may deduce that:
\begin{align*}
|f_*(j_n(S_n))|& \leq |f|\left(\frac{|j_n(S_n)|}{(s_{n+1/2}-s_{n+1})} \right) \leq  2\frac{|j_n(S_n)|^2}{(s_{n+1/2}-s_{n+1})^2}\\
                  & \leq   2e^{2p}\frac{\|j_n\|^2|S_n|^2}{(s_n-s_{n+1/2})^{2p}(s_{n+1/2}-s_{n+1})^2}\\ 
\end{align*}
The sequence $a'$ with terms
$$  a'_n:=2e^{2p}\frac{\|j_n\|^2}{(s_n-s_{n+1/2})^{2p}(s_{n+1/2}-s_{n+1})^2}$$
is subquadratic with order $\a$. Therefore the right-hand side of our estimate is of the form $a_n'|S_n|^2$ with $a'=(a_n') \in \SM^+(\a)$.

 The power-series:
$$e^{-z} \in \CM\{z\}$$
is the Borel transform of 
\[ \frac{1}{1+z} \in \CM\{z\} , \]
therefore 
the remainder term of the iteration satisfies the estimate:
\begin{align*}
|e^{-[j_n(S_n),-]}\s_n(S_n) |& <  2\frac{\|j_n\|\, \| \s_n\|\,  |S_n|}{(s_{n+1/2}-s_{n+1})^{p+1}}.
 \end{align*}
The sequence  with terms
$$ b_n= 2\frac{\|j_n\|\, \| \s_n\|\,  }{(s_{n+1/2}-s_{n+1})^{p+1}}$$
is now subquadratic with order $\b>\a$.
Combining the two estimates we obtain an estimate of the form:
\begin{align*}
|S_{n+1}|  &<  a_n'|S_n|^2+ b_n|S_n|. 
 \end{align*} 
where $a' \in \SM^+$ and $b \in \SM^-$ are subquadratic with orders $\a<\b$.
The sequence with terms
$$a''_n=\frac{s_{n+1/2}-s_{n+1}}{2\|j_n\| }$$
is a falling subquadratic sequence with exponent $\a$. Thus we may choose $\e$ such that for any $n \in \NM$:
$$\e b_n \leq a''_n  .$$
The sequence $(x_k)$ defined by
\begin{align*}
x_0&=|S_0|\\
x_{n+1} &=a_n'  x_n^2+ b_n x_n
\end{align*}
majorates $|S_n|$ and,  according to \ref{SSS::model}, there exists $\dt$ such that for $x_0<\dt$ the sequence is bounded by $\e b$. In particular the iteration $(S_n)$ is well defined and smaller than $\e b$. This concludes the proof of the 
convergence theorem.
\end{proof}
%%%%%%%%%%%%%%%%%%%%%%%%%%%%%%%%%%%%%%%%%%%%%%%%%%
\subsection{Proof of theorem on the stability of quasi-periodic motions}
 Like in the formal case, we start from a perturbed  quasi-periodic motion 
 $$X=X_{\nu}+S_0=\sum_{i=1}^{2d}\nu_i(\dt)\theta_i+S_0$$
and consider the corresponding perturbation of the versal unfolding 
\[X_0=\Vt+S_0.\] 
We choose a subquadratic sequence $a \leq \s(\nu(0))$ defining a set $W(a)$. By Theorem~\ref{T::normal_form} and \ref{SSS::exponential} then for any $\e$ there exists $\dt$ such that
for $|S_0| \leq \dt $ the sequence 
 \[ e^{-[v_n,-]}e^{-[v_{n-1},-]} \dots e^{-[v_0,-]} \]
converges to a mapping 
$$\psi  \in {\Ht}om(\Theta^k(W),\Theta^k(W) ),\ k \geq 1$$
such that
$$\psi(\Vt+S_0)=\Vt.$$
and $\psi $ is $\e$-close to the identity in the $C^k$-topology.
Choose Whitney extensions $r_1,\dots,r_{2d}$ of the functions
$$R_i(\dt,\phi):=\psi(\phi_i),\ i=1,\dots,2d $$
and let $A$ be the Jacobian matrix $(\d_{\phi_j} R_i)$ at the origin.
For $\e$ small enough this matrix is invertible and the Jacobian of
$\phi-A^{-1}r $ at the origin $r=(r_1,\dots,r_n)$ can be made arbitrarily small.  By the implicit function theorem, we may find
a $C^k$-map $g=(g_1,\dots,g_{2d})$,
defined in a neighbourhood
of the origin such that
$$r(\dt,\phi)=0 \iff \phi=g(\dt) .$$
As we remarked previously, the image of $X_\nu+S_0$ under $\psi$ is the restriction of $\Vt$ to
$$R_1=\dots=R_{2d}=0. $$ 
This means that
$$(*)\quad \psi(X_\nu+S_0)= \sum_{i=1}^{2d}(\nu_i(\dt)+g_i(\dt))\theta_i$$
By assumption the map
$$\dt \mapsto \nu(\dt) $$
is a submersion at the origin and for $\e$ small enough
the Jacobian matrix 
$$D\nu(0)+Dg(0)=A^{-1} (\d_{\dt_j} r_i)$$
has maximal rank at the origin and therefore $\nu +g$ is a submersion in a sufficiently small neighbourhood of the origin.
By Ehresmann's lemma, the map $\nu$ defines, over the set $Z(a)_{\infty,s}$,
a local fibration by invariant tori and these carry a quasi-periodic motion by $(*)$.
 The sets of $Z(a)_{\infty,s}$ are of positive measure therefore these invariant tori form a positive measure set as well. This proves the theorem.
\bibliographystyle{plain}
\bibliography{torus_submitted_Duke.bib}

\begin{thebibliography}{10}

\bibitem{Arnold_KAM}
V.I. Arnold.
\newblock { Proof of a theorem of A. N. Kolmogorov on the preservation of
  conditionally periodic motions under a small perturbation of the
  Hamiltonian}.
\newblock {\em Uspehi Mat. Nauk}, 18(5):13--40, 1963.
\newblock English translation: Russian Math. Surveys.

\bibitem{Broer_Huitema_Sevryuk_book}
H.~W. Broer, G.B. Huitema, and M.B. Sevryuk.
\newblock {\em Quasi-periodic motions in families of dynamical systems: order
  amidst chaos}.
\newblock Springer, 2009.

\bibitem{Broer_Huitema_Takens}
H.~W. Broer, G.B Huitema, F.~Takens, and B.L.J. Braaksma.
\newblock {\em Unfoldings and bifurcations of quasi-periodic tori}.
\newblock Number 421. American Mathematical Soc., 1990.

\bibitem{Broer_Huitema_Sevryuk_families}
H.W. Broer, G.B. Huitema, and M.B. Sevryuk.
\newblock Families of quasi-periodic motions in dynamical systems depending on
  parameters.
\newblock In {\em Nonlinear Dynamical Systems and Chaos}, pages 171--211.
  Springer, 1996.

\bibitem{Bruno}
A.D. Bruno.
\newblock {Analytic form of differential equations I}.
\newblock {\em Trans. Moscow Math. Soc.}, 25:131--288, 1971.

\bibitem{Douady_these}
A.~Douady.
\newblock Le probl\`eme des modules pour les sous-espaces analytiques compacts
  d'un espace analytique donn\'e.
\newblock {\em Annales de l'Institut Fourier}, 16(1):1--95, 1966.

\bibitem{Fejoz_KAM}
J.~F{\'e}joz.
\newblock {A proof of the invariant torus theorem of Kolmogorov}.
\newblock {\em Regular and Chaotic Dynamics}, 17(1):1--5, 2012.

\bibitem{arithmetic}
M.~Garay.
\newblock {Arithmetic Density}.
\newblock {\em Proceedings of the Edinburgh Mathematical Society},
  59(3):691--700, 2016.

\bibitem{KAM_book}
M.~Garay and D.~van Straten.
\newblock {KAM Theory}.
\newblock ArXiv: 1204.2493, 2018.

\bibitem{Gonzalez_Haro_DelaLLave}
A.~Gonz{\'a}lez-Enr{\'\i}quez, A.~Haro, and R.~De~la Llave.
\newblock {\em {Singularity theory for non-twist KAM tori}}, volume 227 of {\em
  Memoirs of the AMS}.
\newblock American Mathematical Society, 2014.

\bibitem{Hauser_these}
H.~Hauser.
\newblock Sur la construction de la d\'eformation semi-universelle d'un germe
  d'espace analytique.
\newblock {\em Ann. \'Ec. Norm. Sup. Paris}, 18:1--56, 1985.

\bibitem{Herman_KAM}
M.~Herman.
\newblock Exemples de flots {H}amiltoniens dont aucune perturbation en
  topologie {$C^{\infty}$} n'a d'orbites p\'eriodiques sur un ouvert de
  surfaces d'\'energies.
\newblock {\em C. R. Acad. Sci. Paris. S\'erie I}, \textbf{312}(13):989--994,
  1991.

\bibitem{Kolmogorov_KAM}
A.N. Kolmogorov.
\newblock {On the conservation of quasi-periodic motions for a small
  perturbation of the Hamiltonian function}.
\newblock {\em Dokl. Akad. Nauk SSSR}, 98:527--530, 1954.

\bibitem{Moser_convergent}
J.~Moser.
\newblock Convergent series expansions for quasi-periodic motions.
\newblock {\em Math. Annalen}, 169(1):136--176, 1967.

\bibitem{Moser_KAM}
J.~Moser.
\newblock {On the construction of almost periodic solutions for ordinary
  differential equations (Tokyo, 1969)}.
\newblock In {\em Proc. Internat. Conf. on Functional Analysis and Related
  Topics}, pages 60--67. Univ. of Tokyo Press, 1969.

\bibitem{Poschel_lecture}
J.~P\"oschel.
\newblock {A lecture on the classical KAM theorem}.
\newblock In A.B. Katok, R.~de~la Llave, Ya.B. Pesin, and H.Weiss, editors,
  {\em Smooth Ergodic Theory and Its Applications}, volume~69 of {\em
  Proceedings of Symposia in Pure Mathematics}, pages 707--732. Amer. Math.
  Soc., 2001.

\bibitem{Russmann_degeneracy}
H.~R\"ussmann.
\newblock {Invariant tori in non-degenerate nearly integrable Hamiltonian
  systems}.
\newblock {\em Regular and chaotic dynamics}, 6(2):119--204, 2001.

\bibitem{Sevryuk_KAM}
M.B. Sevryuk.
\newblock {The classical KAM theory at the dawn of the twenty-first century}.
\newblock {\em Moscow Math. Journal}, 3(3):1113--1144, 2003.

\bibitem{Yoccoz_Herman}
J.-C. Yoccoz.
\newblock {Travaux de Herman sur les tores invariants}.
\newblock {\em Ast\'erique}, 206:311--344, 1992.

\end{thebibliography}
\end{document}